\documentclass[hidelinks,onefignum,onetabnum]{siamart220329}
\usepackage{lipsum}
\usepackage{amsfonts}
\usepackage{graphicx}
\usepackage{epstopdf}
\usepackage{algorithmic}
\ifpdf
  \DeclareGraphicsExtensions{.eps,.pdf,.png,.jpg}
\else
  \DeclareGraphicsExtensions{.eps}
\fi

\newsiamremark{remark}{Remark}
\newsiamremark{hypothesis}{Hypothesis}
\crefname{hypothesis}{Hypothesis}{Hypotheses}
\newsiamremark{assumption}{Assumption}
\crefname{assumption}{Assumption}{Assumption}
\newsiamthm{claim}{Claim}

\headers{Stability and convergence of MPC}{D. W. M. Veldman and E. Zuazua}

\title{Local Stability and Convergence of Unconstrained Model Predictive Control\thanks{Submitted to the editors on \today.
\funding{This project has received funding from the European Research Council (ERC) under the European Union’s Horizon 2020 research and innovation programme (grant agreement No. 694126-DyCon and the Marie Sklodowska-Curie grant agreement No. 765579-ConFlex),  the Alexander von Humboldt-Professorship program, the Transregio 154 Project ‘‘Mathematical Modelling, Simulation and Optimization Using the Example of Gas Networks’’, project C08, of the German DFG, the grant PID2020-112617GB-C22, ‘‘Kinetic equations and learning control’’ of the Spanish MINECO, and the COST Action grant CA18232, ‘‘Mathematical models for interacting dynamics on networks’’ (MAT-DYN-NET).}}}

\author{Dani\"el W. M. Veldman\thanks{Chair for Dynamics, Control and Numerics, Alexander von Humboldt Professorship, Friedrich-Alexander-Universit\"at Erlangen-N\"urnberg, Cauerstrasse 11, 91058, Erlangen, Germany 
  (\email{daniel.wm.veldman@fau.de}, \email{enrique.zuazua@fau.de}).}
\and Enrique Zuazua\footnotemark[2] \thanks{Chair in Computational Mathematics, Fundaci\'on Deusto, Av. de las Universidades 24, 48007, Bilbao, Spain and Departamento de Matem\'aticas, Universidad Autonoma de Madrid, 28049, Madrid, Spain.}}

\usepackage{amsopn}

\usepackage{amsmath}
\usepackage{amssymb}
\usepackage{subcaption}

\ifpdf
\hypersetup{
  pdftitle={Local Stability and Convergence of Unconstrained Model Predictive Control},
  pdfauthor={Daniel W.M. Veldman and Enrique Zuazua}
}
\fi

\begin{document}

\maketitle

% REQUIRED
\begin{abstract}
The local stability and convergence for Model Predictive Control (MPC) of unconstrained nonlinear dynamics based on a linear time-invariant plant model is studied. Based on the long-time behavior of the solution of the Riccati Differential Equation (RDE), explicit error estimates are derived that clearly demonstrate the influence of the two critical parameters in MPC: the prediction horizon $T$ and the control horizon $\tau$. In particular, if the MPC-controller has access to an exact (linear) plant model, the MPC-controls and the corresponding optimal state trajectories converge exponentially to the solution of an infinite-horizon optimal control problem when $T-\tau \rightarrow \infty$. When the difference between the linear model and the nonlinear plant is sufficiently small in a neighborhood of the origin, the MPC strategy is locally stabilizing and the influence of modeling errors can be reduced by choosing the control horizon $\tau$ smaller. The obtained convergence rates are validated in numerical simulations.
\end{abstract}

% REQUIRED
\begin{keywords}
Convergence, Model Predictive Control, Receding Horizon Control, Stability
\end{keywords}

% REQUIRED
\begin{MSCcodes}
49N10, 93D15
\end{MSCcodes}

\section{Introduction}
Model Predictive Control (MPC) is a well-established and widely-used feedback control strategy which has received a vast amount of attention in the last four decades, see for example the survey papers \cite{garcia1989,allgower1999,mayne2000,lee2011} and the books \cite{grune2017, rawlings2019}. The main advantages of the MPC paradigm are that 1) the feedback nature of MPC creates, unlike classical optimal control theory (see, e.g., \cite{lee1967}), additional robustness against disturbances, modeling errors, and implementation errors and that 2) MPC can, unlike many other techniques for feedback control design (see, e.g., \cite{skogestad2005}), be applied to nonlinear systems with input and state constraints. Additionally, MPC reduces the horizon over which optimal control problems need to be solved, and may therefore reduce the memory requirements and computational cost for the implementation of the controller. 

The idea for MPC can already be found in the classical book by Lee and Markus \cite{lee1967}:
\begin{quote}
One technique for obtaining a feedback controller synthesis from knowledge of open-loop controllers is to measure the current control process state and then compute very rapidly for the open-loop control function. The first portion of this function is then used during a short time interval, after which a new measurement of the process state is made and a new open-loop control function is computed for this new measurement. The procedure is then repeated. 
\end{quote}
Due to the limited available computational power at that time, it was difficult to compute the open-loop control function `very rapidly' and this observation did not receive much attention initially. The continual increase in computational power has enabled the application of MPC in industrial applications since the late 1970's and MPC has received a huge interest since then, both from industry and academia. 

Two central questions in the literature on MPC are whether the MPC-feedback is stabilizing and how the performance of the MPC-controller compares to the optimal performance. Classically, the stability question has been addressed by imposing proper terminal constraints or terminal costs, see, e.g., \cite{allgower1999}, but later research has shown that this is in fact not necessary. In particular, Lars Gr\"une \cite{grune2009} proposed an analysis method for MPC that does not require state constraints or terminal costs to guarantee stability or performance. These ideas have been particularly influential in the last decade. The original paper \cite{grune2009} only considers discrete-time systems, but the ideas have been extended to continuous-time systems in \cite{reble2012}. A peculiar artifact in the estimates from \cite{reble2012} is that they blow up when the control horizon approaches zero. As is also remarked in \cite{reble2012}, this behavior is counterintuitive and has stimulated some research in MPC with short control horizons (also called instant-MPC), see, e.g., \cite{yoshida2019}. 

In this paper, an analysis method for MPC that is based on the (exponential) convergence of the solution of the Riccati Differential Equation (RDE) to the (symmetric positive-definite) solution of the Algebraic Riccati Equation (ARE) is proposed. In contrast to many existing results, the presented analysis only requires standard controllability and observability assumptions and no terminal constraints or terminal costs. In contrast to the existing results in \cite{reble2012}, the estimates in this paper also remain bounded (and actually improve) when the control horizon approaches zero. 

The remainder of this paper is structured as follows. In Section \ref{sec:main}, the MPC strategy is introduced and the main ideas and results of this paper are summarized. Section \ref{sec:proofs} contains the detailed proofs of the results from Section \ref{sec:main}. Section \ref{sec:numerics} contains two numerical examples that validate the results in Section \ref{sec:main}. Finally, conclusions and discussions are presented in Section \ref{sec:conc}. 
The discrete-time analogues of the developments in this paper can be found in Appendix \ref{sec:discretetime}.

\section{Main ideas and results} \label{sec:main}

\subsection{Model predictive control} \label{ssec:MPC}
Consider the nonlinear dynamical system
\begin{equation}
\dot{y}(t) = f(y(t),u(t)) + w(t), \qquad y(0) = y_0, \label{eq:dyn_y}
\end{equation} 
where the state $y(t)$ evolves in $\mathbb{R}^n$ starting from the initial condition $y_0$, the control $u(t)$ evolves in $\mathbb{R}^m$ (with $m \leq n$), the disturbance $w(t)$ evolves in $\mathbb{R}^n$, and $f : \mathbb{R}^n \times \mathbb{R}^m \rightarrow \mathbb{R}^n$ is Lipschitz and satisfies $f(0,0) = 0$. In many practical situations, the dynamical system \eqref{eq:dyn_y} is not available for control but the state $y(t)$ can be measured at certain time instances $t = k\tau$ (for some $\tau > 0$ and $k \in \mathbb{N}$). 
The goal in MPC is therefore to find a control $u(t)$ that stabilizes the system based on these measurements and an imperfect plant model. In this paper, it is assumed that the plant model available for control is linear and time-invariant (LTI)
\begin{equation}
\dot{x}(t) = Ax(t) + Bu(t), \label{eq:dyn_x}
\end{equation}
where the state $x(t)$ also evolves in $\mathbb{R}^n$, $A \in \mathbb{R}^{n \times n}$ is the system matrix, and $B \in \mathbb{R}^{n \times m}$ is the input matrix. Ideally, $A = f_y(0,0)$  and $B = f_u(0,0)$ but these conditions may be violated because $f$ is typically not known exactly. The goal in MPC is to use the measurements $y(k\tau)$ and the linear model \eqref{eq:dyn_x} find a control $u(t)$ that locally stabilizes the origin of the nonlinear plant \eqref{eq:dyn_y}, preferably with a nearly minimal cost
\begin{equation}
I_\infty(u) = \frac{1}{2}\int_0^\infty \left( |Cy(t)|^2 + (u(t))^\top R u(t) \right) \ \mathrm{d}t, \label{eq:def_I}
\end{equation}
where $C \in \mathbb{R}^{\ell \times n}$ is the output matrix (with $\ell \leq n$), $R \in \mathbb{R}^{m\times m}$ is a symmetric positive definite matrix, and $y(t)$ satisfies \eqref{eq:dyn_y}. 
Our results indicate that MPC can achieve this goal locally when $f_y(0,0) - A$, $f_u(0,0) - B$, and $w(t)$ are sufficiently small.

\begin{remark} \label{rem:RHCvsMPC}
Receding Horizon Control (RHC) is closely related to MPC and the terms are often used interchangeably in the literature. The convention from \cite{ito2002, azmi2016, azmi2018} is also adopted in this paper: the term RHC is reserved for MPC based on a perfect plant model. In particular, RHC refers to MPC with $f(y,u) = Ay + Bu$ and $w \equiv 0$. 
\end{remark}

To introduce the MPC algorithm, fix a prediction horizon $T \geq \tau$ and introduce $u_T^*(\cdot;x_1,t_1) \in L^2([t_1, t_1+T], \mathbb{R}^m)$ as the minimizer of the functional $J_T(u;x_1,t_1)$
\begin{equation}
J_T(u) = 
\frac{1}{2}(x(t_1+T))^\top E_T x(t_1+T) + \frac{1}{2} \int_{t_1}^{t_1+T} \left( |C x(t)|^2 + (u(t))^\top R u(t) \right) \ \mathrm{d}t, \label{eq:JT}
\end{equation}
subject to the dynamics \eqref{eq:dyn_x} and $x(t_1) = x_1$. 
Here, $C$ and $R$ are as in \eqref{eq:def_I} and the terminal cost $E_T \in \mathbb{R}^{n\times n}$ is a symmetric positive semi-definite matrix that can improve the stability and convergence of MPC, see, e.g., \cite{allgower1999, ito2002} and Remark \ref{rem:K0ET} below. The dependence of quantities in the optimal control problem \eqref{eq:JT} on $x_1$ and $t_1$ will be omitted throughout this paper when no confusion can occur. 

Furthermore, let $x^*_T(t;x_1,t_1)$ and $y^*_T(t;x_1,t_1)$ denote the trajectories of \eqref{eq:dyn_x} and \eqref{eq:dyn_y} resulting from the control $u^*_T(t;x_1,t_1)$
\begin{align}
\dot{x}^*_T(t) &= Ax^*_T(t) + Bu^*_T(t), \quad & & x^*_T(t_1) = x_1, \label{eq:dyn_xstar} \\
\dot{y}^*_T(t) &= f(y^*_T(t), u^*_T(t)) + w(t), & & y^*_T(t_1) = x_1. \label{eq:dyn_ystar}
\end{align}
In RHC,  $f(y,u) = Ay + Bu$ and $w(t) \equiv 0$ so that $x^*_T=y^*_T$, see Remark \ref{rem:RHCvsMPC}.

The MPC control $u_{\mathrm{MPC}}(t)$ is now constructed as follows. Measure $y_{\mathrm{MPC}}(0) = y_0$, compute the control $u^*_T(t;y_{\mathrm{MPC}}(0),0)$ on $[0,T]$, set $u_{\mathrm{MPC}}(t) = u^*_T(t;y_{\mathrm{MPC}}(0),0)$ on $(0,\tau]$ and apply it to the plant \eqref{eq:dyn_y}. This leads to the state trajectory $y_{\mathrm{MPC}}(t) = y^*_T(t;0,y_0)$ on $[0,\tau]$. Next, measure $y_{\mathrm{MPC}}(\tau)$, compute $u^*_T(t;y_{\mathrm{MPC}}(\tau), \tau)$ on $[\tau, \tau+T]$, set $u_{\mathrm{MPC}}(t) = u^*_T(t;y_{\mathrm{MPC}}(\tau), \tau)$ on $(\tau, 2\tau]$, and apply it to the plant \eqref{eq:dyn_y}. This leads to the trajectory $y_{\mathrm{MPC}}(t) = y^*_T(t; y_{\mathrm{MPC}}(\tau),\tau)$ on $[\tau, 2\tau]$ which can again be measured at time $t = 2\tau$. Repeating this procedure results in Algorithm \ref{alg:MPC}. 

\begin{algorithm}
\textbf{Step 1} Choose a prediction horizon $T < \infty$ and a control horizon $0 < \tau \leq T$. Set $k = 0$. 

\textbf{Step 2} Measure $y_{\mathrm{MPC}}(k\tau)$ and compute the control $u^*_T(t; y_{\mathrm{MPC}}(k\tau), k\tau)$.

\textbf{Step 3} Apply the control $u_{\mathrm{MPC}}(t) = u^*_T(t; y_{\mathrm{MPC}}(k\tau), k\tau)$ to the plant \eqref{eq:dyn_y} during $t \in (k\tau, k\tau+\tau]$. 

\textbf{Step 4} Increase $k$ by 1 and go to step 2. 

\caption{Model Predictive Control}
\label{alg:MPC}
\end{algorithm}

\begin{figure}
\centering
\includegraphics[width=0.7\columnwidth]{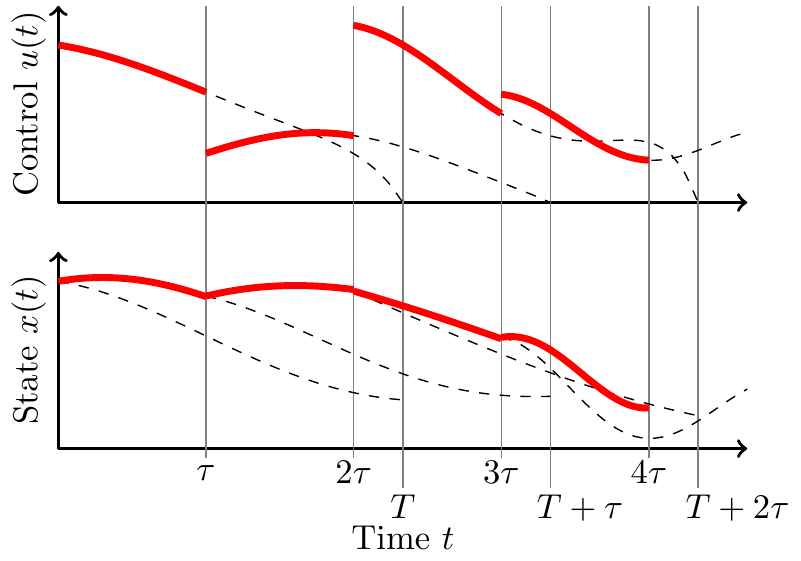}
\caption{The MPC strategy in Algorithm \ref{alg:MPC}. The black dashed lines indicate $u_T^*(t; y_{\mathrm{MPC}}(k\tau), k\tau)$ and $x^*_T(t; y_{\mathrm{MPC}}(k\tau),k\tau)$ on the time interval $[k\tau, k\tau + T]$. The red thick lines show $u_T^*(t; y_{\mathrm{MPC}}(k\tau),k\tau)$ and $y_T^*(t; y_{\mathrm{MPC}}(k\tau),k\tau)$ in the time interval $[k\tau, k\tau+\tau)$, which together form $u_{\mathrm{MPC}}(t)$ and $y_{\mathrm{MPC}}(t)$.}
\label{fig:MPC}
\end{figure}

Figure \ref{fig:MPC} shows a typical control $u_{\mathrm{MPC}}(t)$ and state trajectory $y_{\mathrm{MPC}}(t)$ resulting from Algorithm \ref{alg:MPC}. Note that $y_{\mathrm{MPC}}(t)$ is continuous, but that the control $u_{\mathrm{MPC}}(t)$ is not. The controls $u^*_T(t; y_{\mathrm{MPC}}(k\tau), k\tau)$ in Figure \ref{fig:MPC} (the dashed lines in top graph) vanish at $t = k\tau + T$. This represents the situation in which the terminal cost $E_T = 0$. Furthermore, note that the trajectories $x^*_T(t;y_{\mathrm{MPC}}(k\tau), k\tau)$ (dashed lines in bottom graph) differ from $y_{\mathrm{MPC}}(t)$ (red line in bottom graph). 

The central question in the paper is now for which prediction horizons $T$ and control horizons $\tau$ Algorithm \ref{alg:MPC} is stabilizing and how this answer depends on the modeling errors $f(y,u) - Ay - Bu$ and $w(t)$. 

\subsection{Main idea} \label{ssec:Riccati}
The analysis in this paper is based on the observation that in RHC (i.e., MPC with $f(y,u) = Ay + Bu$ and $w(t) \equiv 0$), the control $u_{\mathrm{RHC}}(t)$ and the corresponding state trajectory $y_{\mathrm{RHC}}(t)$ can be considered as approximations of the control $u^*_\infty(t)$ and corresponding state trajectory $x^*_\infty(t)$ that minimize 
\begin{equation}
J_\infty(u) = \frac{1}{2}\int_{0}^\infty \left( |Cx(t)|^2  +  (u(t))^\top R u(t) \right) \ \mathrm{d}t, \label{eq:Jinf}
\end{equation}
subject to \eqref{eq:dyn_x} and $x(0) = y_0$.
 If $(A,B)$ is controllable and $(A,C)$ is observable, it is well-known that the optimal trajectory $x^*_\infty(t)$  is given by, see, e.g., \cite{sage1968}
\begin{equation}
\dot{x}_\infty^*(t) = A_\infty x_\infty^*(t), \qquad x_\infty^*(t_0) = y_0, \label{eq:xinf_dyn}
\end{equation}
where $A_\infty$ is
\begin{equation}
A_\infty := A - BR^{-1}B^\top P_\infty, \label{eq:Ainfty}
\end{equation}
with $P_\infty$ the unique symmetric positive-definite solution of the Algebraic Riccati Equation (ARE)
\begin{equation}
A^\top P_\infty + P_\infty A - P_\infty B R^{-1} B^\top P_\infty + C^\top C = 0. \label{eq:ARE}
\end{equation}
Because the controllability of $(A,C)$ implies that the matrix $A_\infty$ in \eqref{eq:Ainfty} is Hurwitz, see, e.g., \cite[Lemma 2.6]{porretta2013}, there exist a growth bound $\mu_\infty > 0$ and overshoot constant $M_\infty \geq 1$ such that for all $t \geq 0$
\begin{equation}
\| e^{A_\infty t}\| \leq M_\infty e^{-\mu_\infty t}, \label{eq:growthbound}
\end{equation}
where $\| \cdot \|$ denotes the operator norm. In particular, $x^*_\infty(t) \rightarrow 0$ for $t \rightarrow \infty$. 

It is also well-known (see, e.g., \cite{sage1968}) that the finite-horizon optimal control problem is solvable and that the corresponding state trajectory $x^*_T(t,x_1,t_1)$ satisfies
\begin{equation}
\dot{x}^*_T(t) = (A - B R^{-1}B^\top P_T(t-t_1))x^*_T(t), \qquad\qquad x^*_T(t_1) = x_1,
\label{eq:xT_dyn}
\end{equation}
where $P_T(t)$ is the $\mathbb{R}^{n\times n}$-matrix valued solution of the Ricatti Differential Equation (RDE) (on $[0,T]$)
\begin{equation}
- \dot{P}_T(t) = A^\top P_T(t) + P_T(t) A - P_T(t) B R^{-1} B^\top P_T(t) + C^\top C, \quad P_T(T) = E_T.  \label{eq:RDE}
\end{equation}
Note that \eqref{eq:RDE} is solved backward in time starting from the final condition.
It is therefore convenient to introduce $\mathcal{P}(t)$ (for all $t \geq 0$) as the solution of 
\begin{equation}
\dot{\mathcal{P}}(t) = A^\top \mathcal{P}(t) + \mathcal{P}(t) A
- \mathcal{P}(t) B R^{-1} B^\top \mathcal{P}(t) + C^\top C, \qquad \mathcal{P}(0) = E_T.  \label{eq:RDE_reversed}
\end{equation}
Comparing \eqref{eq:RDE} and \eqref{eq:RDE_reversed}, it follows that
\begin{equation}
P_T(t) = \mathcal{P}(T-t). \label{eq:P_def}
\end{equation}

The main idea for the analysis of this paper is now the following. In RHC, $y_{\mathrm{RHC}}(t) = y^*_T(t;y_{\mathrm{RHC}}(k\tau),k\tau) = x^*_T(t;y_{\mathrm{RHC}}(k\tau),k\tau)$ for $t \in [k\tau, k\tau+\tau]$. Because $x^*_T(t; x_1, t_1)$ satisfies \eqref{eq:xT_dyn}, it follows that in RHC
\begin{equation}
\dot{y}_{\mathrm{RHC}}(t) = A_{T,\tau}(t) y_{\mathrm{RHC}}(t), \qquad y_{\mathrm{RHC}}(0) = y_0, \label{eq:xRH_ATtau}
\end{equation}
where the $\tau$-periodic matrix $A_{T,\tau}(t)$ is defined as
\begin{equation}
A_{T,\tau}(t) := A - BR^{-1}B^\top P_{T,\tau}(t), \label{eq:def_ATt}
\end{equation}
with $P_{T,\tau}(t)$ being the $\tau$-periodic matrix
\begin{equation}
P_{T,\tau}(t) = P_T(t \mod \tau) = \mathcal{P}(T - (t \mod \tau)). \label{eq:def_PTt} 
\end{equation}

The following lemma shows that $\mathcal{P}(t) \rightarrow P_\infty$ for $t \rightarrow \infty$ and is fundamental for the analysis in this paper. 
\begin{lemma} \label{lem:Pconv}
Assume that $(A,B)$ is controllable and $(A,C)$ is observable. 
Let $P_\infty$ be the symmetric positive definite solution of the ARE \eqref{eq:ARE}, let $\mathcal{P}(t)$ be the solution of the time-reversed RDE \eqref{eq:RDE_reversed}, and let $\mu_\infty > 0$ be the growthbound in \eqref{eq:growthbound}, then there exists a constant $K_0$ (independent of $t$) such that
\begin{equation}
\| \mathcal{P}(t) - P_\infty \| \leq K_0 e^{-2\mu_\infty t}. \label{eq:Pconv}
\end{equation}
\end{lemma}
A brief proof is given in Appendix \ref{app:RDE}, but similar results can be found in \cite{callier1994, porretta2013}. 

Note that $\mathcal{P}(t) \rightarrow P_\infty$ for $t \rightarrow \infty$ clearly implies that $A_{T,\tau}(t) \rightarrow A_\infty$ for $T - \tau \rightarrow \infty$ (see \eqref{eq:def_ATt} and \eqref{eq:Ainfty}). Because $A_\infty$ is Hurwitz, $A_{T,\tau}(t)$ will be Hurwitz for all time $t$ if $T - \tau$ is sufficiently large. This observation is the key to establish the stability and convergence results for RHC and MPC in the next subsection. 

\begin{remark} \label{rem:K0ET}
The proof in Appendix \ref{app:RDE} shows that $K_0 \geq M_\infty^2 \| P_\infty - E_T \|$ and that $K_0 = 0$ when $E_T = P_\infty$. %The latter follows because $\min_u J_\infty(u; x_1, 0) = x_1^\top P_\infty x_1$. Bellman's optimality principle then shows that the minimizer of $J_T(u; y_0, 0)$ in \eqref{eq:JT} is equal to the minimizer of $J_\infty(u)$ (restricted to $[0,T]$) and the MPC control is equal to the optimal control in this case. In general, $E_T$ should thus be chosen such that $x_1^\top E_T x_1$ is a good approximation of $\min_u J_\infty(u; x_1, 0)$ (for all $x_1$). 
\end{remark}

\begin{remark} \label{rem:stab_detect}
The assumptions in Lemma \ref{lem:Pconv} can be relaxed to $(A,B)$ being stabilizable and $(A,C)$ being detectable, see \cite{callier1994}. 
\end{remark}

\subsection{Main Results} \label{ssec:main}
The results in this section are based on the exponential convergence of the solution to the RDE to the solution of the ARE from Lemma \ref{lem:Pconv} and the observation that $y_{\mathrm{RHC}}(t)$ is described by the periodic feedback law in \eqref{eq:xRH_ATtau}. All results are based on the following assumption that enables us to use Lemma \ref{lem:Pconv}. 
\begin{assumption}
The pair $(A,B)$ is controllable and the pair $(A,C)$ is observable. 
\end{assumption}
This assumption can be relaxed slightly, see Remark \ref{rem:stab_detect}. For clarity, proofs are postponed to Section \ref{sec:proofs}.  

\subsubsection{Results for RHC} \label{ssec:RHC}
Our first result is a stability result for RHC. 
\begin{theorem}[Stability of RHC] \label{thm:RHCstab}
If $f(y,u) = Ay + Bu$ and $w \equiv 0$, there exists a constant $K_1$ independent of $t$, $y_0$, $T$, and $\tau$ such that
\begin{equation}
|y_{\mathrm{RHC}}(t)| \leq M_\infty e^{-\mu_{T-\tau} t} |y_0|, \label{eq:thm_RHCstab}
\end{equation}
where
\begin{equation}
\mu_{T-\tau} := \mu_\infty - K_1 e^{-2\mu_\infty(T-\tau)}. \label{eq:def_muTtau}
\end{equation}
\end{theorem}
RHC is thus stabilizing when $\mu_{T-\tau} > 0$, i.e.\  when $T-\tau$ sufficiently large. 

Because Lemma \ref{lem:Pconv} implies that $A_{T,\tau}(t)$ in \eqref{eq:def_ATt} converges to $A_\infty$ in \eqref{eq:Ainfty} for $T-\tau \rightarrow \infty$, $y_{\mathrm{RHC}}(t)$ in \eqref{eq:xRH_ATtau} converges to $x^*_\infty(t)$ in \eqref{eq:xinf_dyn} for $T-\tau \rightarrow \infty$. The corresponding control $u_{\mathrm{RHC}}(t)$ also converges to $u^*_\infty(t)$. These ideas are made precise in the following theorem. 

\begin{theorem}[Convergence of RHC] \label{thm:RHCconv}
If $f(y,u) = Ay + Bu$ and $w \equiv 0$, there exists a constant $K$ independent of $t$, $T$, $y_0$ and $\tau$ such that
\begin{equation}
|y_{\mathrm{RHC}}(t) - x^*_\infty(t) | + |u_{\mathrm{RHC}}(t) - u^*_\infty(t)| \leq K e^{-2\mu_\infty(T-\tau)} \left(1 + t \right)e^{-\mu_{T-\tau}t} |y_0|. \label{eq:thm_uRHC}
\end{equation}
If $\mu_{T-\tau} > 0$, there exists a constant $K$ independent of $t$, $T$, $y_0$, and $\tau$ such that
\begin{equation}
J_\infty(u_{\mathrm{RHC}}) - J_\infty(u^*_\infty) \leq K  \frac{e^{-4\mu_\infty(T-\tau)}}{\mu_{T-\tau}^3} |y_0|^2. \label{eq:thm_JRHC}
\end{equation}
\end{theorem}

\begin{remark} \label{rem:existing}
The suboptimality estimates for RHC from \cite{reble2012, azmi2016, azmi2018} take the form
\begin{equation}
\alpha_{T,\tau} J_\infty (u_{\mathrm{RHC}}) \leq J_\infty(u^*_\infty), \label{eq:suboptimality}
\end{equation}
where $\alpha_{T,\tau} \in (-\infty, 1]$. Note that $\alpha_{T,\tau} = 1$ implies that the performance of the RHC is optimal and that \eqref{eq:suboptimality} does not provide any information when $\alpha_{T,\tau} < 0$.
A typical estimate for $\alpha_{T,\tau}$ is of the form (see \cite[Section 2]{azmi2016})
\begin{equation}
\alpha_{T,\tau} = 1 - \frac{1}{\tau} \frac{1}{T-\tau} \frac{(\gamma(T))^2}{\alpha_\ell^2},
\end{equation}
for a certain bounded function $\gamma(T)$ and a coefficient $\alpha_\ell > 0$ that should satisfy $|Cx|^2 \geq \alpha_\ell|x|^2$ for all $x$. It thus follows that $\alpha_{T,\tau} \rightarrow -\infty$ if either $\tau \rightarrow 0$, or $\tau \rightarrow T$, or $C$ becomes singular. In contrast, the estimates in Theorems \ref{thm:RHCstab} and \ref{thm:RHCconv} do not blow up for $\tau \rightarrow 0$ or $\tau \rightarrow T$ and only require that $(A,C)$ is observable, and are thus also applicable in situations in which $C$ is not invertible.  
\end{remark}

\subsubsection{Results for MPC}
The analysis of the MPC algorithm involves the Lipschitz constant $L$ of $f(y,u) - Ay - Bu$, i.e., for all $(y,u),(y',u') \in \mathbb{R}^n \times \mathbb{R}^m$, 
\begin{equation}
|f(y,u) - Ay - Bu - f(y',u') + Ay' + Bu'| < L(|y-y'| + |u-u'|). \label{eq:Lipschitz}
\end{equation}
Note that $L = 0$ when $f(y,u) = Ay + Bu$ and that $L = \max\{ \| \Delta A \|, \| \Delta B \| \}$ when $f(y,u) = (A + \Delta A)y + (B + \Delta B)u$. 

\begin{remark}
The stability estimates below give conditions on $L$, $T$, and $\tau$ for which $y_{\mathrm{MPC}}(t)$ remains bounded. When these conditions are satisfied, \eqref{eq:Lipschitz} only needs to be satisfied in a neighborhood of the origin and the results below then also only hold for sufficiently small initial conditions $y_0$ and disturbances $w(t)$. Note that if $f$ is $C^2$, $A = f_y(0,0)$, and $B = f_u(0,0)$, \eqref{eq:Lipschitz} can be achieved for any $L > 0$ in a sufficiently small neighborhood of the origin. 
\end{remark} 

With this notation we obtain the following stability result for MPC. 
\begin{theorem}[Stability of  MPC] \label{thm:MPCstab}
There exist constants $K_1$, $K_2$, and $K$ independent of $t$, $y_0$, $T$, $\tau$, $L$, and $w(t)$, such that
\begin{multline}
|y_{\mathrm{MPC}}(t)| \leq M_\infty e^{-\mu_{L,T,\tau} t} |y_0| \\ + \frac{1-e^{-\mu_{L,T,\tau} t}}{\mu_{L,T,\tau}} K (1 + (L+1)\tau e^{K(L+1)\tau})|w|_{L^\infty(0,t)}, \label{eq:thm_MPCstab}
\end{multline}
where
\begin{equation}
\mu_{L,T,\tau} := \mu_\infty - K_1 e^{-2\mu_\infty(T-\tau)} - K_2 L - KL(L+1)\tau e^{K(L+1)\tau} . \label{eq:def_muL}
\end{equation}

\end{theorem}

Note that $\mu_{L,T,\tau}$ does not depend on $w(t)$ and that $\mu_{L,T,\tau} = \mu_{T-\tau}$ in \eqref{eq:def_muTtau} when $L = 0$, i.e.\ when $f(y,u) = Ay + Bu$. Note that the closed-loop MPC dynamics is Input-to-State Stable (ISS) w.r.t.\ the disturbance $w$ if $\mu_{L,T,\tau} > 0$. If $K_2 L < \mu_\infty$, $\mu_{L,T,\tau}$ is positive for $T-\tau$ sufficiently large and $\tau$ sufficiently small. MPC can thus only be stabilizing when the modeling errors (measured by $L$) are small enough.

The following convergence result shows that, for $T-\tau \rightarrow \infty$ and $\tau \rightarrow 0$, $y_{\mathrm{MPC}}(t)$ converges to the solution $z_\infty(t)$ of
\begin{equation}
\dot{z}_\infty(t) = f(z_\infty(t), -R^{-1}B^\top P_\infty z_\infty(t)) +  w(t), \qquad\qquad z_\infty(0) = y_0, \label{eq:dyn_zinf}
\end{equation}
and the corresponding control $u_{\mathrm{MPC}}(t)$ converges to 
\begin{equation}
v_\infty(t) := -R^{-1}B^\top P_\infty z_\infty(t). \label{eq:vinf}
\end{equation}
Note that $z_\infty(t)$ is the trajectory resulting the application of infinite-horizon feedback operator for the linear model \eqref{eq:dyn_x} to the nonlinear model \eqref{eq:dyn_y}. 
Note in particular that $v_\infty(t)$ is not the minimizer of $I_\infty$ in \eqref{eq:def_I}, see Remark \ref{rem:vinfty} below. 

\begin{theorem}[Convergence of MPC] \label{thm:MPCconv}
If $K_2 L < \mu_\infty$, there exist a constant $K$ independent of $t$, $y_0$, $T$, $\tau$, $L$, and $w(t)$ such that
\begin{multline}
|y_{\mathrm{MPC}}(t) - z_\infty(t)| + |u_{\mathrm{MPC}}(t) - v_\infty(t)| \\
\leq Ke^{-2\mu_\infty(T-\tau)} \left( \frac{L+1}{\mu_\infty - K_2 L}|y_{\mathrm{MPC}}|_{L^1(0,t)} + |y_{\mathrm{MPC}}(t)|\right)\\
+ K \tau e^{K(L+1)\tau} \frac{L+1}{\mu_\infty - K_2 L} \left(L |y_{\mathrm{MPC}}|_{L^1(0,t)} + |w|_{L^\infty(0,t)}  \right). \label{eq:MPCconv}
\end{multline}
\end{theorem}

Note that Theorem \ref{thm:MPCstab} shows that $|y_{\mathrm{MPC}}|_{L^1(0,t)}$ and $|y_{\mathrm{MPC}}(t)|$ can be bounded for $T-\tau$ sufficiently large and $\tau$ sufficiently small because $K_2L < \mu_\infty$. 

\begin{remark} \label{rem:vinfty}
To understand why $u_{\mathrm{MPC}}(t)$ does not converge to the minimizer of $I_\infty$ in \eqref{eq:def_I}, consider the situation where $f(y,u) = Ay + Bu$ and the disturbance $w(t)$ is nonzero. The minimizer of the infinite-horizon problem should be introduced carefully because there might not exist a control that makes the infinite horizon cost finite. However, by considering problems on a finite horizon $T$ and taking the limit $T \rightarrow \infty$, it can be shown that the optimal state trajectories converge to the solution of (see e.g.\ \cite[Section 5.2]{sage1968})
\begin{equation}
\dot{y}_\infty^*(t) = (A - BR^{-1}B^\top P_\infty) y_\infty^*(t) - BR^{-1}B^\top \xi(t)+ w(t), \qquad y^*_\infty(0) = x_0,  \label{eq:dyn_yinf}
\end{equation}
where $\xi(t)$ is the solution of
\begin{equation}
-\dot{\xi}(t) = (A - BR^{-1}B^\top P_\infty)^\top \xi(t) - P_\infty w(t), \qquad\qquad \xi(\infty) = 0. \label{eq:dyn_xi}
\end{equation}
Note that $\xi(t)$ depends on $w(s)$ for $s \geq t$, about which the MPC-controller has no information at time $t$. It is therefore not possible that $y_{\mathrm{MPC}}(t)$ converges to $y^*_\infty(t)$. 
\end{remark}

\section{Proofs} \label{sec:proofs} This section contains the proofs of Theorems \ref{thm:RHCstab}, \ref{thm:RHCconv}, \ref{thm:MPCstab}, and \ref{thm:MPCconv}. 

\subsection{Stability of RHC (Theorem \ref{thm:RHCstab})} \label{ssec:RHCstab}

Lemma \ref{lem:Pconv} and \eqref{eq:def_PTt} show that
\begin{equation}
\max_{t\geq 0} \| P_{T,\tau}(t)  - P_\infty \| \leq K_0 e^{-2\mu_\infty (T-\tau)}. \label{eq:PTdiff}
\end{equation}
From the definitions of $A_{T,\tau}(t)$ in \eqref{eq:def_ATt} and $A_\infty$ in \eqref{eq:Ainfty}, it thus follows that
\begin{equation}
\max_{t \geq 0} \|A_{T,\tau}(t) - A_\infty \| = \max_{t \geq 0} \| BR^{-1}B^\top (P_{T,\tau}(t) - P_\infty) \| 
\leq K_1' e^{-2\mu_\infty(T-\tau)}, \label{eq:ATtdiff}
\end{equation}
where $K_1' = \| BR^{-1}B^\top\|K_0$. 
Now observe that \eqref{eq:xRH_ATtau} shows that
\begin{equation}
\dot{y}_{\mathrm{RHC}}(t) = A_{T,\tau}(t) y_{\mathrm{RHC}}(t) = A_\infty y_{\mathrm{RHC}}(t) + (A_{T,\tau}(t) - A_\infty) y_{\mathrm{RHC}}(t). 
\end{equation}
The variation of constants formula thus shows that
\begin{equation}
y_{\mathrm{RHC}}(t) = e^{A_\infty t} y_0 + \int_0^t e^{A_\infty(t-s)} (A_{T,\tau}(s) - A_\infty) y_{\mathrm{RHC}}(s) \ \mathrm{d}s. 
\end{equation}
Using the triangle inequality, \eqref{eq:growthbound}, and \eqref{eq:ATtdiff}, it follows that
\begin{equation}
|y_{\mathrm{RHC}}(t)| \leq M_\infty e^{-\mu_\infty t}|y_0| +   \int_0^t e^{-\mu_\infty(t-s)} K_1 e^{-2\mu_\infty(T-\tau)} |y_{\mathrm{RHC}}(s)| \ \mathrm{d}s, \label{eq:xRHbound}
\end{equation}
where $K_1 = M_\infty K_1'$. 
Multiplying \eqref{eq:xRHbound} by $e^{\mu_\infty t}$ and writing $\hat{y}_{\mathrm{RHC}}(t) = e^{\mu_\infty t} y_{\mathrm{RHC}}(t)$, it follows that
\begin{equation}
|\hat{y}_{\mathrm{RHC}}(t)| \leq M_\infty |y_0| + K_1 e^{-2\mu_\infty(T-\tau)} \int_0^t |\hat{y}_{\mathrm{RHC}}(s)| \ \mathrm{d}s. 
\end{equation}
Gr\"onwall's lemma then yields
\begin{equation}
|\hat{y}_{\mathrm{RHC}}(t)| \leq M_\infty |y_0|e^{t K_1e^{-2\mu_\infty(T-\tau)}}. 
\end{equation}
Theorem \ref{thm:RHCstab} now follows after noting that $|y_{\mathrm{RHC}}(t)| = e^{-\mu_\infty t} |\hat{y}_{\mathrm{RHC}}(t)|$. 

\subsection{Convergence of RHC (Theorem \ref{thm:RHCconv})} Throughout the proof, $K$ denotes a generic constant that does not depend on $t$, $y_0$, $T$, and $\tau$ that may vary from line to line. 
Denote $e_{\mathrm{RHC}}(t) := y_{\mathrm{RHC}}(t) - x_\infty^*(t)$ and observe that
\begin{equation}
\dot{e}_{\mathrm{RHC}}(t) = A_{T,\tau}(t) y_{\mathrm{RHC}}(t) - A_\infty x^*_\infty(t) = A_\infty e_{\mathrm{RHC}}(t) + (A_{T,\tau}(t) - A_\infty)y_{\mathrm{RHC}}(t),
\end{equation}
and that $e_{\mathrm{RHC}}(0) = 0$. Therefore,
\begin{equation}
e_{\mathrm{RHC}}(t) = \int_0^t e^{A_\infty(t-s)} (A_{T,\tau}(s) - A_\infty)y_{\mathrm{RHC}}(s) \ \mathrm{d}s. \label{eq:eRH}
\end{equation}
Taking norms, making use of \eqref{eq:growthbound} and \eqref{eq:ATtdiff}, it follows that
\begin{align}
|e_{\mathrm{RHC}}(t)| &\leq K e^{-2\mu_\infty(T-\tau)} \int_0^t e^{-\mu_\infty(t-s)} |y_{\mathrm{RHC}}(s)| \ \mathrm{d}s \nonumber \\
&\leq K e^{-2\mu_\infty(T-\tau)} \int_0^t e^{-\mu_\infty(t-s)} e^{-\mu_{T-\tau}s} |y_0| \ \mathrm{d}s \leq K e^{-2\mu_\infty(T-\tau)} t e^{-\mu_{T-\tau}t} |y_0|, \label{eq:eRHnorm}
\end{align}
where the second inequality follows from the stability result in Theorem \ref{thm:RHCstab} and the third inequality because $\mu_\infty \geq \mu_{T-\tau}$. 

For the bound on the controls, note that comparing \eqref{eq:dyn_x} and \eqref{eq:xinf_dyn} yields $u^*_\infty(t) = -R^{-1}B^\top P_\infty x^*_\infty(t)$. Similarly, comparing \eqref{eq:dyn_xstar} and \eqref{eq:xT_dyn} noting that $x^*_T(t) = y^*_T(t)$ in RHC, yields $u_{\mathrm{RHC}}(t) = -R^{-1}B^\top P_{T,\tau}(t) y_{\mathrm{RHC}}(t)$. Therefore,
\begin{align}
|u_{\mathrm{RHC}}(t) - & u^*_\infty(t)|
= |R^{-1}B^\top (P_{T,\tau}(t)y_{\mathrm{RHC}}(t) - P_\infty x^*_\infty(t))| \nonumber \\
&= |R^{-1}B^\top ((P_{T,\tau}(t) - P_\infty) y_{\mathrm{RHC}}(t) + P_\infty e_{\mathrm{RHC}}(t))| \nonumber \\
&\leq  \| R^{-1}B^\top \| \|P_{T,\tau}(t) - P_\infty \| |y_{\mathrm{RHC}}(t)| + \|R^{-1}B^\top P_\infty  \| |e_{\mathrm{RHC}}(t)| \nonumber \\
&\leq K e^{-2\mu_\infty(T-\tau)} e^{-\mu_{T-\tau}t} |y_0| + K e^{-2\mu_\infty(T-\tau)} t e^{-\mu_{T-\tau}t} |y_0|,   \label{eq:vRH} 
\end{align}
where the last inequality follows after using \eqref{eq:PTdiff} and Theorem \ref{thm:RHCstab} for the first term and \eqref{eq:eRHnorm} for the second term. 

For \eqref{eq:thm_JRHC}, note that because $I_\infty$ is quadratic and $u^*_\infty$ is the minimizer of $I_\infty$
\begin{align}
I_\infty(u_{\mathrm{RHC}}) - I_\infty(u^*_\infty) &= \delta I_\infty(u_\infty^*;v_{\mathrm{RHC}}) + \tfrac{1}{2} \delta^2 I_\infty(u^*_\infty, v_{\mathrm{RHC}}, v_{\mathrm{RHC}}) \nonumber \\
&= \frac{1}{2}\int_0^\infty \left( |C e_{\mathrm{RHC}}(t)|^2 + (v_{\mathrm{RHC}}(t))^\top R v_{\mathrm{RHC}}(t) \right) \ \mathrm{d}t \nonumber \\ 
&\leq \tfrac{1}{2}\| C \|^2 |e_{\mathrm{RHC}}|_{L^2(0,\infty)}^2 + \tfrac{1}{2}\| R \| |v_{\mathrm{RHC}}|_{L^2(0,\infty)}^2. \label{eq:thmJRH_step1}
\end{align}
where $v_{\mathrm{RHC}}(t) = u_{\mathrm{RHC}}(t) - u^*_\infty(t)$, $\delta I_\infty(u;v)$ denotes the Fr\'echet derivative of $I_\infty$ in the point $u$ in the direction $v$, $\delta^2 I_\infty(u;v,v)$ is the Hessian in the point $u$ in the direction $v$. Note that the second equality follows because $\delta I_\infty(u^*_\infty, v) = 0$ for all $v$ (because $u^*_\infty$ is the minimizer of $I_\infty$) and the form of $I_\infty$ in \eqref{eq:def_I}. 

Bounds for $|e_{\mathrm{RHC}}|_{L^2(0,\infty)}^2$ and $|v_{\mathrm{RHC}}|_{L^2(0,\infty)}^2$ follow after squaring \eqref{eq:eRHnorm} and \eqref{eq:vRH} and integrating from $t = 0$ to $t = \infty$ using that
\begin{align}
\int_0^\infty t^2 e^{-2\mu_{T-\tau}t} \ \mathrm{d}t &= \frac{1}{4 \mu_{T-\tau}^3}, \\
\int_0^\infty t e^{-2\mu_{T-\tau}t} \ \mathrm{d}t &= \frac{1}{4 \mu_{T-\tau}^2} < \frac{\mu_\infty}{4 \mu_{T-\tau}^3}, \\
\int_0^\infty e^{-2\mu_{T-\tau}t} \ \mathrm{d}t &= \frac{1}{2 \mu_{T-\tau}} < \frac{\mu_\infty^2}{2 \mu_{T-\tau}^3},
\end{align}
where it is used that $0 < \mu_{T-\tau} < \mu_\infty$. 

\subsection{Stability of MPC (Theorem \ref{thm:MPCstab})} \label{ssec:MPCstab} Throughout the proof, $K$ denotes a generic constant that does not depend on $t$, $y_0$, $T$, $\tau$, $L$, and $w(t)$ that may vary from line to line. 
Write $\tau_t = t - (t \mod \tau)$ and note that Algorithm \ref{alg:MPC} and \eqref{eq:dyn_ystar} yield
\begin{align}
\dot{y}_{\mathrm{MPC}}(t) &= f(y_{\mathrm{MPC}}(t),u_T^*(t)) + w(t) \nonumber \\
&= A_\infty y_{\mathrm{MPC}}(t) + (A_{T,\tau}(t) - A_\infty) y_{\mathrm{MPC}}(t) - BR^{-1}B^\top P_{T,\tau}(t) \varepsilon(t) \nonumber \\
&\qquad + f(y_{\mathrm{MPC}}(t),u_T^*(t)) - Ay_{\mathrm{MPC}}(t) -Bu^*_T(t) + w(t), \label{eq:dyn_yMPC2}
\end{align}
with $A_{T,\tau}(t)$ as in \eqref{eq:def_ATt}, $u^*_T(t) = -R^{-1}B^\top P_{T,\tau}(t) x^*_T(t)$, and $\varepsilon(t) = x^*_T(t) - y_{\mathrm{MPC}}(t)$ with $x_T^*(t) = x_T^*(t, y_{\mathrm{MPC}}(\tau_t), \tau_t)$. Note that \eqref{eq:Lipschitz} (with $(y',u') = (0,0)$) shows that
\begin{align}
|f(y_{\mathrm{MPC}}(t),u^*_T(t)) - A y_{\mathrm{MPC}}(t) &-Bu^*_T(t) | \leq L (|y_{\mathrm{MPC}}(t)| + |u^*_T(t)| ) \nonumber \\
&\leq L ((1+K_2')|y_{\mathrm{MPC}}(t)| + K_2'|\varepsilon(t)| ), \label{eq:Lipschitz_applied}
\end{align}
where $K_2'= \| R^{-1} B^\top \|(K_0 + \| P_\infty \|)$. Here, the second inequality uses that $u^*_T(t) = -R^{-1}B^\top P_{T,\tau}(t) (y_{\mathrm{MPC}}(t) + \varepsilon(t))$ and that $\| P_{T,\tau}(t) \| \leq K_0 + \| P_\infty \|$ by Lemma \ref{lem:Pconv}. Applying the variation of constants formula to \eqref{eq:dyn_yMPC2} using that $y_{\mathrm{MPC}}(0) = y_0$ and taking norms using \eqref{eq:growthbound}, \eqref{eq:ATtdiff}, and \eqref{eq:Lipschitz_applied}, it follows that
\begin{align}
|y_{\mathrm{MPC}}(t)| \leq & M_\infty e^{-\mu_\infty t} |y_0| + \left(K_1 e^{-2\mu_\infty (T-\tau)} + K_2 L \right) \int_0^t e^{-\mu_\infty(t-s)} |y_{\mathrm{MPC}}(s)|  \ \mathrm{d}s \nonumber \\ 
& + K(L+1) \int_0^t e^{-\mu_\infty(t-s)} |\varepsilon(s)| \ \mathrm{d}s + M_\infty \int_0^t e^{-\mu_\infty(t-s)} |w(s)| \ \mathrm{d}s, \label{eq:yMPC_est1}
\end{align}
where $K_2 = M_\infty(1+K_2')$ and $K_1 = M_\infty \| BR^{-1}B^\top \| K_0$ as in Subsection \ref{ssec:RHCstab}.

To estimate $\varepsilon(t)$, note that subtracting \eqref{eq:dyn_yMPC2} from \eqref{eq:xT_dyn} shows that
\begin{align}
\dot{\varepsilon}(t) &= A_\infty\varepsilon(t) + (A_{T,\tau}(t) - A_\infty)\varepsilon(t) - BR^{-1}B^\top P_{T,\tau}(t) \varepsilon(t) \nonumber \\ 
&\qquad - f(y_{\mathrm{MPC}}(t),u^*_T(t))+ A y_{\mathrm{MPC}}(t) + B u^*_T(t) - w(t), \qquad\varepsilon(\tau_t) = 0. \label{eq:thm_MPCstab_step97}
\end{align}
Applying the variation of constants formula and taking norms yields
\begin{align}
|\varepsilon(t)| &\leq M_\infty \int_{\tau_t}^t \left| (A_{T,\tau}(s) - A_\infty)\varepsilon(s) - BR^{-1}B^\top P_{T,\tau}(s) \varepsilon(s) \right| \ \mathrm{d}s \nonumber \\ 
&\qquad + M_\infty \int_{\tau_t}^t \left| f(y_{\mathrm{MPC}}(t),u^*_T(t))- A y_{\mathrm{MPC}}(t) - B u^*_T(t) + w(t) \right| \ \mathrm{d}s \\
&\leq K(L+1) \int_{\tau_t}^t  |\varepsilon(s)| \ \mathrm{d}s \nonumber + KL |y_{\mathrm{MPC}}|_{L^1(\tau_t,t)} + \tau M_\infty |w|_{L^\infty(0,t)}. 
\end{align}
where \eqref{eq:growthbound}, \eqref{eq:ATtdiff}, and \eqref{eq:Lipschitz_applied} were used. Gr\"onwall's lemma thus shows that
\begin{equation}
|\varepsilon(t)| \leq e^{K(L+1) \tau} \left( K L |y_{\mathrm{MPC}}|_{L^1(\tau_t,t)} + \tau M_\infty |w|_{L^\infty(0,t)} \right). \label{eq:boundf}
\end{equation}
Next, define $\hat{y}_{\mathrm{MPC}}(t) = e^{\mu_\infty t} y_{\mathrm{MPC}}(t)$ and multiply \eqref{eq:yMPC_est1} by $e^{\mu_\infty t}$ to find
\begin{align}
|\hat{y}_{\mathrm{MPC}}(t)| &\leq M_\infty |y_0| + \left(K_1 e^{-2\mu_\infty (T-\tau)} + K_2L \right)|\hat{y}_{\mathrm{MPC}}|_{L^1(0,t)}  \nonumber \\ 
&\qquad + K (L+1) e^{K(L+1) \tau}  \int_0^t  \left( e^{\mu_\infty \tau} L |\hat{y}_{\mathrm{MPC}}|_{L^1(\tau_s,s)} + e^{\mu_\infty s}\tau|w|_{L^\infty(0,t)} \right)   \ \mathrm{d}s \nonumber \\
& \qquad\qquad\qquad\qquad + M_\infty \int_0^t e^{\mu_\infty s} |w(s)| \ \mathrm{d}s, \label{eq:thm_MPCstab_step94}
\end{align}
where it has been used that
\begin{equation}
e^{\mu_\infty s} |y_{\mathrm{MPC}}|_{L^1(\tau_s,s)} = \int_{\tau_s}^s e^{\mu_\infty (s-r)} |\hat{y}_{\mathrm{MPC}}(r)| \ \mathrm{d}r \leq e^{\mu_\infty \tau} |\hat{y}_{\mathrm{MPC}}|_{L^1(\tau_s,s)}.
\end{equation}
For the integral of $|\hat{y}_{\mathrm{MPC}}|_{L^1(\tau_s,s)}$ in \eqref{eq:thm_MPCstab_step94}, observe that for $t \in [k\tau, k\tau+\tau)$
\begin{align}
&\int_0^t |\hat{y}_{\mathrm{MPC}}|_{L^1(\tau_s,s)} \ \mathrm{d}s = 
\int_{k\tau}^t |\hat{y}_{\mathrm{MPC}}|_{L^1 (k\tau,s)} \ \mathrm{d}s + \sum_{\ell=0}^{k-1} \int_{\ell\tau}^{\ell\tau+\tau} |\hat{y}_{\mathrm{MPC}}|_{L^1(\ell\tau,s)} \ \mathrm{d}s  \nonumber \\
&\qquad\qquad \leq \tau |\hat{y}_{\mathrm{MPC}}|_{L^1(k\tau,t)} + \tau \sum_{\ell=0}^{k-1} |\hat{y}_{\mathrm{MPC}}|_{L^1(\ell\tau,\ell\tau+\tau)} = \tau |\hat{y}_{\mathrm{MPC}}|_{L^1(0,t)}. \label{eq:yMPCbound}
\end{align}
The result then follows by inserting \eqref{eq:yMPCbound} into \eqref{eq:thm_MPCstab_step94}, finding a bound for $|\hat{y}_{\mathrm{MPC}}(t)|$ using  Gr\"onwall's lemma, and using that $|y_{\mathrm{MPC}}(t)| = e^{-\mu_\infty t} |\hat{y}_{\mathrm{MPC}}(t)|$. 

\subsection{Convergence of MPC (Theorem \ref{thm:MPCconv})} Just as in the proof of Theorem \ref{thm:MPCstab}, $K$ denotes a generic constant that does not depend on $t$, $y_0$, $T$, $\tau$, $L$, and $w(t)$ that may vary from line to line. Note that \eqref{eq:dyn_zinf} can be rewritten as
\begin{equation}
\dot{z}_\infty(t) = f(z_\infty(t), v_\infty(t)) - Az_\infty(t) - Bv_\infty(t) + A_\infty z_\infty(t) + w(t),
\end{equation}
where \eqref{eq:Ainfty} and \eqref{eq:vinf} have been used. Writing $e_{\mathrm{MPC}}(t) = y_{\mathrm{MPC}}(t) - z_\infty(t)$ and subtracting this equation from \eqref{eq:dyn_yMPC2}, it follows that
\begin{multline}
\dot{e}_{\mathrm{MPC}}(t) = A_\infty e_{\mathrm{MPC}}(t) + (A_{T,\tau}(t) - A_\infty)y_{\mathrm{MPC}}(t) - BR^{-1}B^\top P_{T,\tau}(t)\varepsilon(t) \\
+  f(y_{\mathrm{MPC}}(t),u^*_T(t)) - f(z_\infty(t), v_\infty(t))- Ae_{\mathrm{MPC}}(t) - B(u^*_T(t) - v_\infty(t)), \label{eq:dyn_eMPC}
\end{multline}
where the notation is the same as in \eqref{eq:dyn_yMPC2}. Note that \eqref{eq:Lipschitz} shows that
\begin{align}
|f(y_{\mathrm{MPC}}(t),u^*_T(t)) & - f(z_\infty(t), v_\infty(t))- Ae_{\mathrm{MPC}}(t) - B(u^*_T(t) - v_\infty(t))| \nonumber \\
&\leq L(|e_{\mathrm{MPC}}(t)| + |u^*_T(t) - v_\infty(t)|) \nonumber \\
&\leq L ((1+K_2') |e_{\mathrm{MPC}}(t)| + K_2'|\varepsilon(t)| + K e^{-2\mu_\infty(T-\tau)} |y_{\mathrm{MPC}}(t)|), \label{eq:Lipschitz_applied2}
\end{align}
where $K'_2 = \| R^{-1}B^\top\|(K_0 + \|P_\infty \|)$ as in Subsection \ref{ssec:MPCstab} and it was used that
\begin{align}
u^*_T(t) - v_\infty(t) &= -R^{-1}B^\top( P_{T,\tau}(t)x^*_T(t) - P_\infty z_\infty(t)) \nonumber \\
&= -R^{-1}B^\top \left( P_{T,\tau}(t) \varepsilon(t) + (P_{T,\tau}(t) - P_\infty)y_{\mathrm{MPC}}(t) + P_\infty e_{\mathrm{MPC}}(t) \right), \label{eq:error_MPCcontrols}
\end{align}
and \eqref{eq:PTdiff} was used to bound $P_{T,\tau}(t) - P_\infty$. 
Applying the variation of constants formula to \eqref{eq:dyn_eMPC} and taking norms using \eqref{eq:growthbound}, \eqref{eq:ATtdiff}, and \eqref{eq:Lipschitz_applied2}, it follows that
\begin{align}
|e_{\mathrm{MPC}}(t)| &\leq K_2 L \int_0^t e^{-\mu_\infty (t-s)} |e_{\mathrm{MPC}}(s)|  \ \mathrm{d}s + K(L+1)e^{-2\mu_\infty(T-\tau)} |y_{\mathrm{MPC}}|_{L^1(0,t)} \nonumber  \\
&\qquad + K(L+1) \int_0^t e^{-\mu_\infty(t-s)} |\varepsilon(s)| \ \mathrm{d}s, \label{eq:eMPCbound1}
\end{align}
where $K_2 = M_\infty(1 + K_2')$ as in Subsection \ref{ssec:MPCstab}. 
For the last term, note that \eqref{eq:boundf} and an estimate similar to \eqref{eq:yMPCbound} show that
\begin{align}
\int_0^t  e^{-\mu_\infty(t-s)} &|\varepsilon(s)|  \ \mathrm{d}s \nonumber \\
&\leq e^{K(L+1)\tau} \int_0^t e^{-\mu_\infty(t-s)} (KL|y_{\mathrm{MPC}}|_{L^1(\tau_s,s)} + \tau M_\infty |w|_{L^\infty(0,t)}) \ \mathrm{d}s \nonumber \\
&\leq K \tau e^{K(L+1)\tau} ( L |y_{\mathrm{MPC}}|_{L^1(0,t)} + |w|_{L^\infty(0,t)}). \label{eq:eMPCbound12}
\end{align}
Inserting \eqref{eq:eMPCbound12} into \eqref{eq:eMPCbound1} yields
\begin{equation}
|e_{\mathrm{MPC}}(t)| \leq K_2 L \int_0^t e^{-\mu_\infty (t-s)} |e_{\mathrm{MPC}}(s)| \ \mathrm{d}s  + \alpha(t), \label{eq:eMPCbound2} 
\end{equation}
with
\begin{equation}
\alpha(t) =  
K(L+1) \left( (e^{-2\mu_\infty(T-\tau)}+L \tau e^{K(L+1)\tau} )|y_{\mathrm{MPC}}|_{L^1(0,t)} 
+ \tau e^{K(L+1)\tau}  |w|_{L^\infty(0,t)} \right). \nonumber
\end{equation}
Applying Gr\"onwall's lemma to \eqref{eq:eMPCbound2}, using that $\mu_\infty > K_2 L$ by assumption, gives
\begin{equation}
|e_{\mathrm{MPC}}(t)| \leq \alpha(t) \left(1 + K_2 L \frac{1 - e^{-(\mu_\infty -K_2 L)t}}{\mu_\infty - K_2 L} \right) \leq \alpha(t) \frac{\mu_\infty}{\mu_\infty - K_2 L}. \label{eq:eMPCbound99}
\end{equation}

The estimate for $u_{\mathrm{MPC}}(t) - v_\infty(t)$ follows by taking norms in \eqref{eq:error_MPCcontrols} using \eqref{eq:boundf} for the first term, \eqref{eq:PTdiff} for the second term, and \eqref{eq:eMPCbound99} for the third term.

\section{Numerical examples} \label{sec:numerics}

This section contains two numerical examples that validate the convergence rates from Theorems \ref{thm:RHCconv} and \ref{thm:MPCconv}. 
\subsection{Example 1}
\begin{figure}
\centering
\includegraphics[width=0.8\columnwidth]{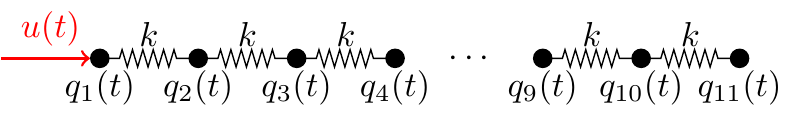}
\caption{The considered system of 11 interconnected point-masses}
\label{fig:model}
\end{figure}

\begin{figure}
\begin{subfigure}{\columnwidth}
\centering
\includegraphics[width=0.75\textwidth]{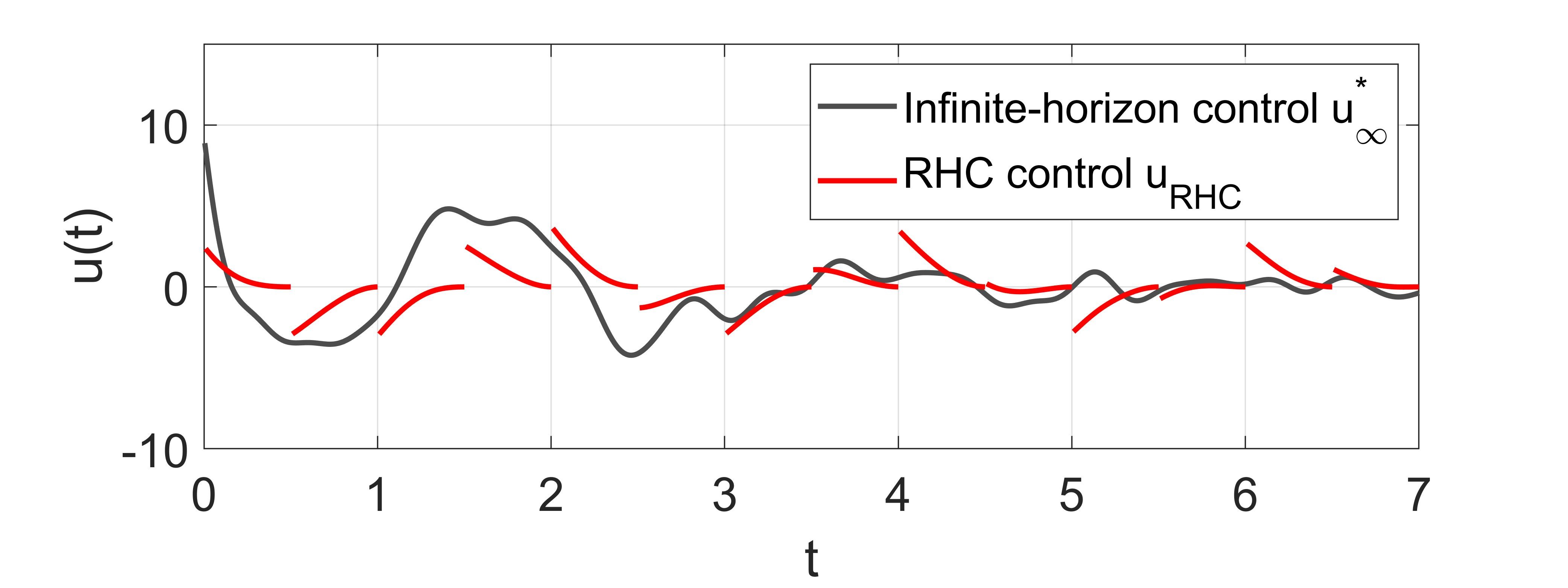}
\caption{$\tau = \tfrac{1}{2}$, $T - \tau = 0$, $f(y,u) = Ay + Bu$, and $w(t) \equiv 0$}
\label{fig:RHC1}
\end{subfigure}
\begin{subfigure}{\columnwidth}
\centering
\includegraphics[width=0.75\textwidth]{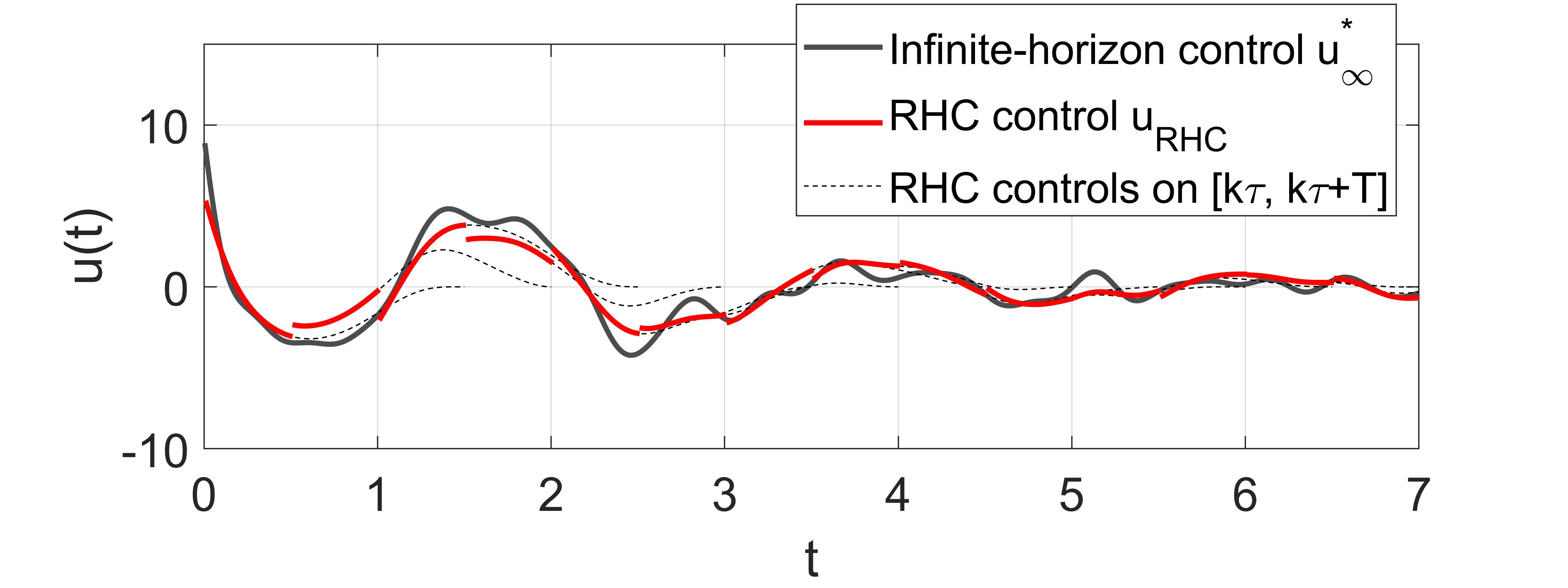}
\caption{$\tau = \tfrac{1}{2}$, $T - \tau = 1$, $f(y,u) = Ay + Bu$, and $w(t) \equiv 0$}
\label{fig:RHC2}
\end{subfigure}
\begin{subfigure}{\columnwidth}
\centering
\includegraphics[width=0.75\textwidth]{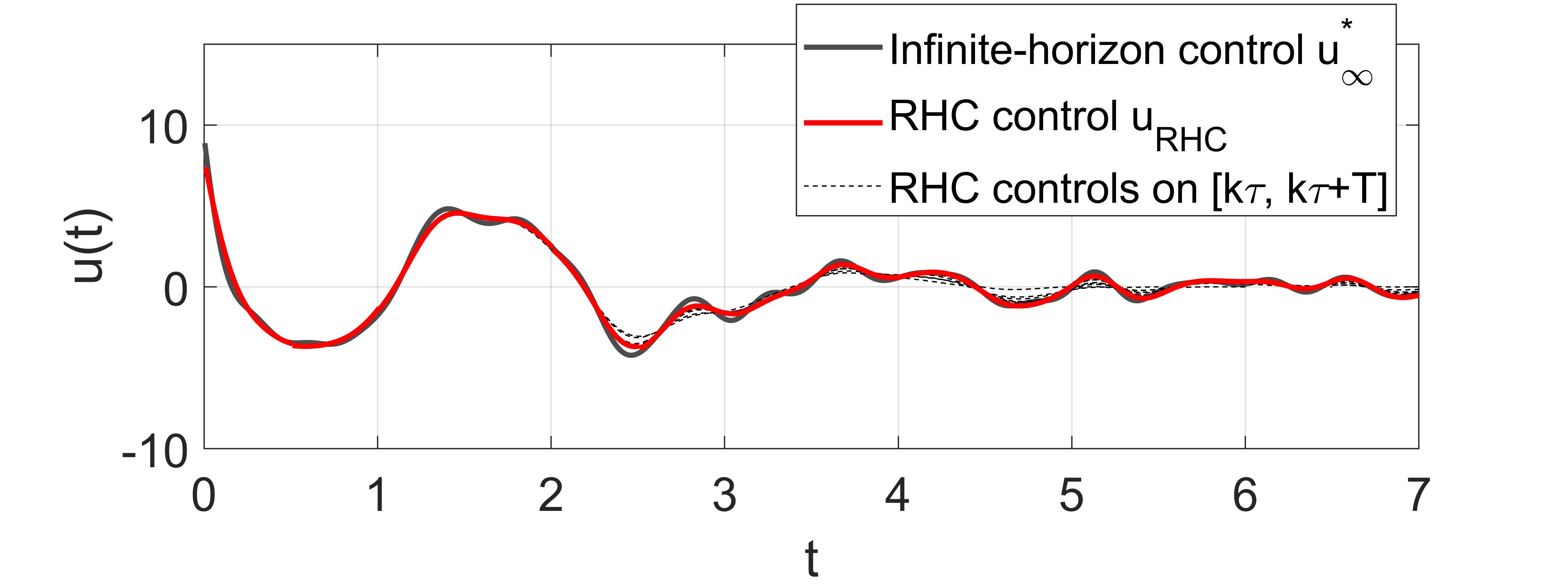}
\caption{$\tau = \tfrac{1}{2}$, $T - \tau = 4$, $f(y,u) = Ay + Bu$, and $w(t) \equiv 0$}
\label{fig:RHC3}
\end{subfigure}
\begin{subfigure}{\columnwidth}
\centering
\includegraphics[width=0.75\textwidth]{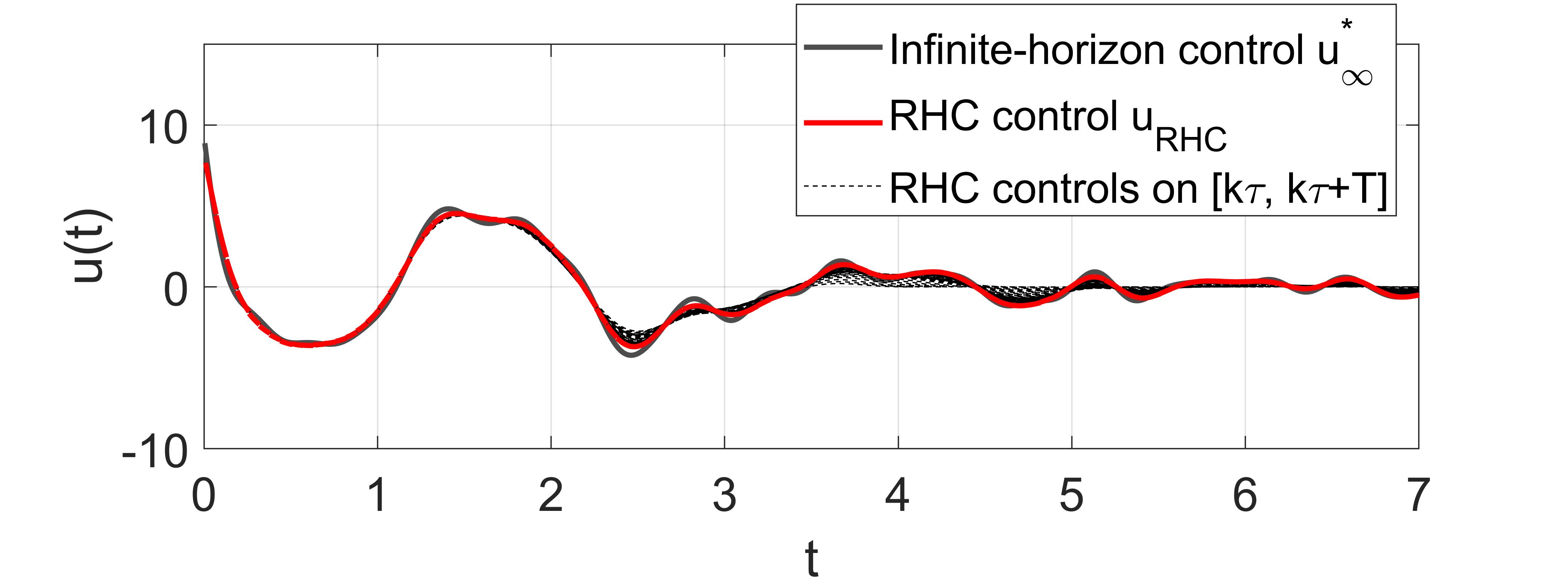}
\caption{$\tau = \tfrac{1}{16}$, $T - \tau = 4$, $f(y,u) = Ay + Bu$, and $w(t) \equiv 0$}
\label{fig:RHC4}
\end{subfigure}
\caption{The controls $u_{\mathrm{RHC}}(t)$ (red) generated by the RHC algorithm for four values of $(T,\tau)$ compared to the infinite-horizon control $u^*_\infty(t)$ (grey). The dashed lines also show the part of controls $u^*_T(t,y_{\mathrm{RHC}}(k\tau), k\tau)$ that are not applied to the plant. }
\label{fig:RHCnum}
\end{figure}

\begin{figure}
\begin{subfigure}{\columnwidth}
\centering
\includegraphics[width=0.75\textwidth]{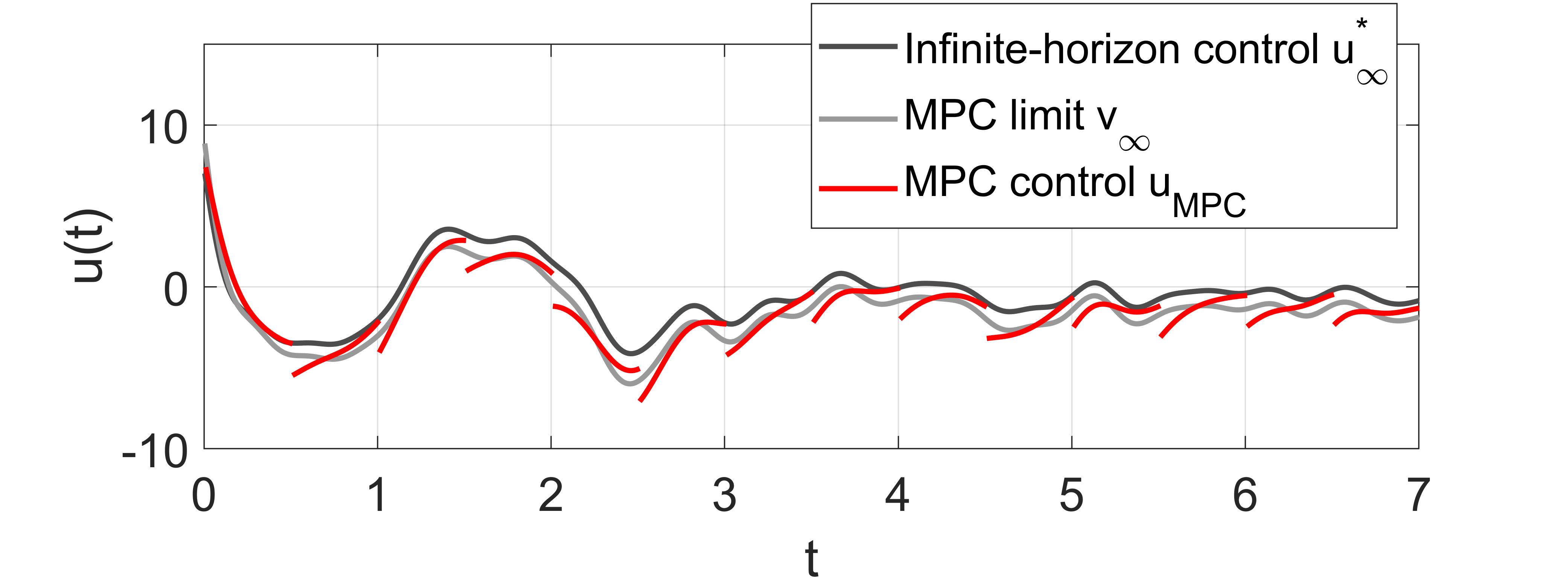}
\caption{$\tau = \tfrac{1}{2}$, $T - \tau = 4$, $f(y,u) = Ay + Bu$, and $w(t) = \bar{w}$}
\label{fig:MPCAw1}
\end{subfigure}
\begin{subfigure}{\columnwidth}
\centering
\includegraphics[width=0.75\textwidth]{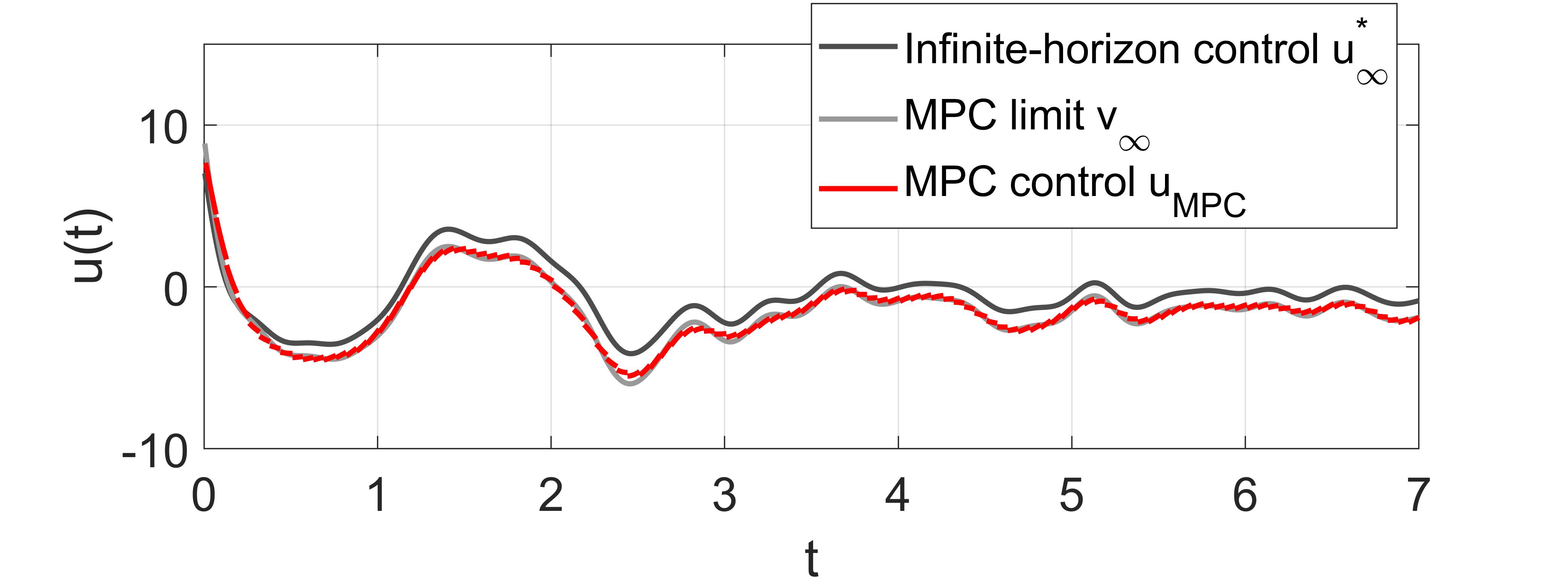}
\caption{$\tau = \tfrac{1}{2}$, $T - \tau = 4$, $f(y,u) = Ay + Bu$, and $w(t) = \bar{w}$}
\label{fig:MPCAw2}
\end{subfigure}

\begin{subfigure}{\columnwidth}
\centering
\includegraphics[width=0.75\textwidth]{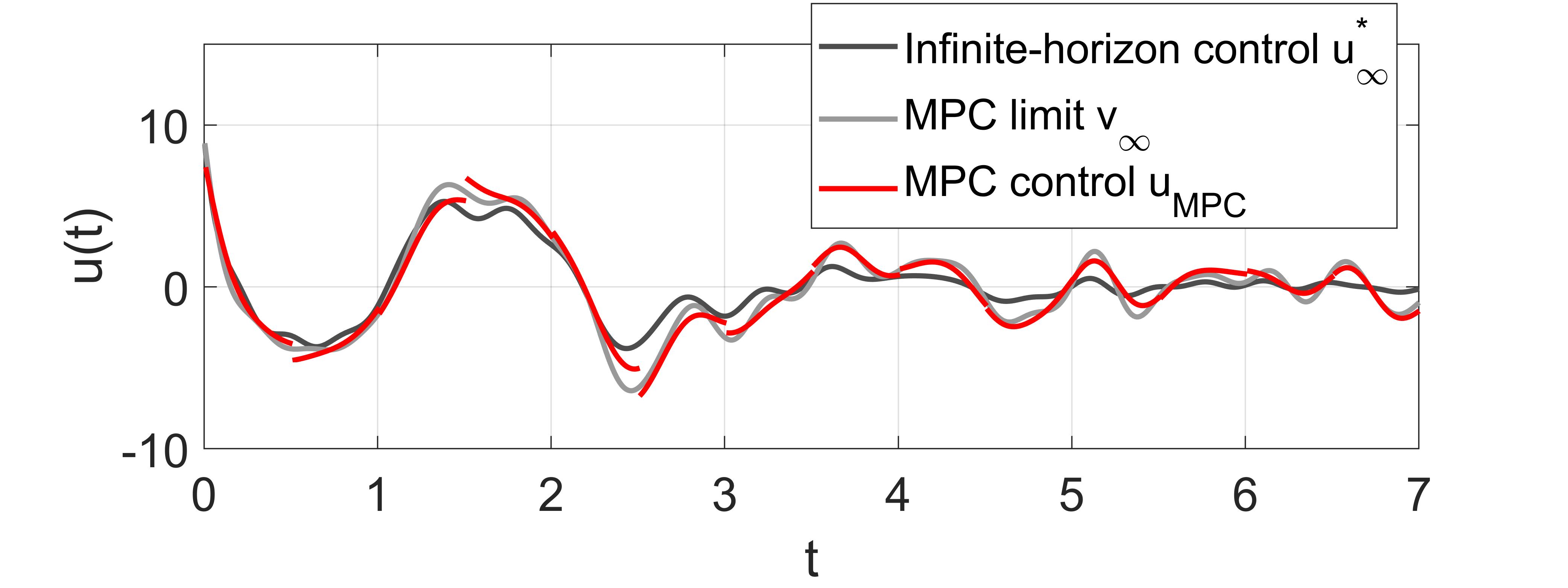}
\caption{$\tau = \tfrac{1}{2}$, $T - \tau = 4$, $f(y,u) = \tilde{A}y + Bu$, and $w(t) \equiv 0$}
\label{fig:MPCAw3}
\end{subfigure}
\begin{subfigure}{\columnwidth}
\centering
\includegraphics[width=0.75\textwidth]{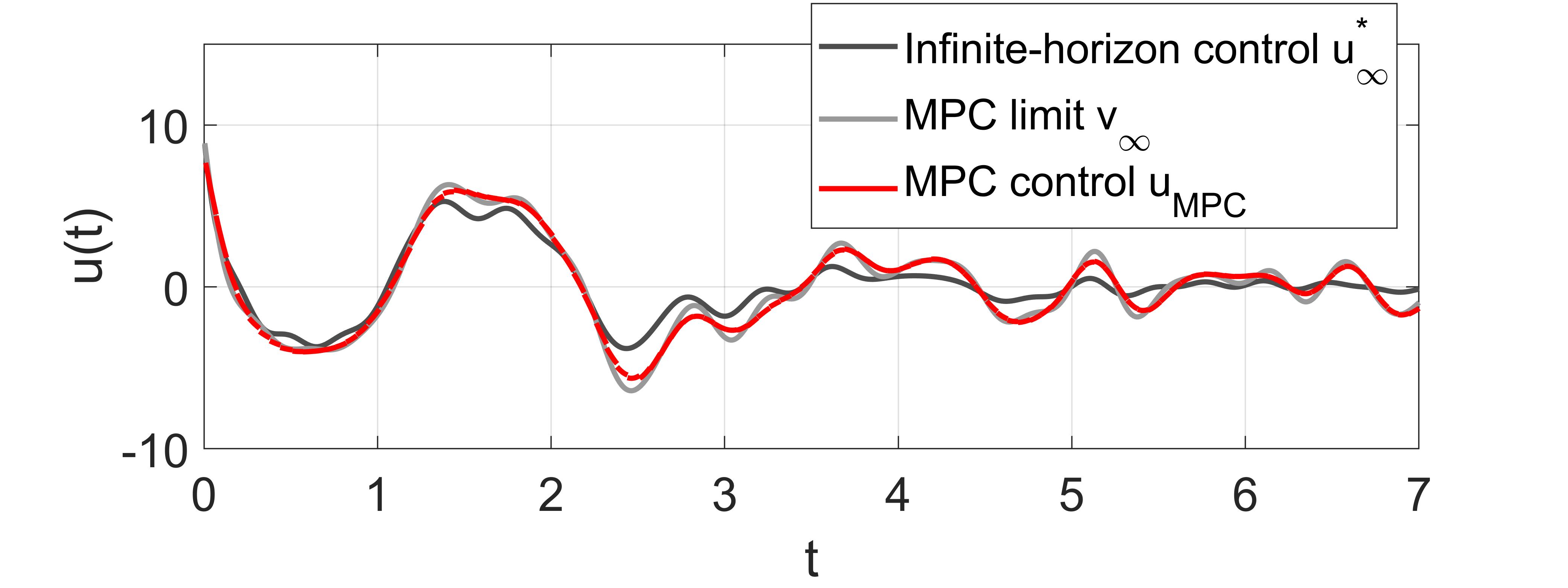}
\caption{$\tau = \tfrac{1}{16}$, $T - \tau = 4$, $f(y,u) = \tilde{A}y + Bu$, and $w(t) \equiv 0$}
\label{fig:MPCAw4}
\end{subfigure}
\caption{The controls $u_{\mathrm{MPC}}(t)$ (red) generated by MPC with a (constant) disturbance, i.e.\  $f(y,u) = Ay + Bu$ and $w(t) = \bar{w}$, and a perturbation in the system matrix, i.e.\ $f(y,u) = \tilde{A}y + Bu$ and $w(t) \equiv 0$, for two values of $(T,\tau)$ compared to the optimal infinite-horizon control $u^*_\infty(t)$ (dark grey) and the MPC limit $v_\infty(t)$ in (light grey).}
\label{fig:MPCAw}
\end{figure}

Consider a system of $N = 11$ point masses of unit mass connected with springs with springconstants $k = 100$ as in Figure \ref{fig:model}. The positions of the point masses (w.r.t.\ an inertial frame) $q_i(t)$ are stored in the vector $q(t)$ of length 11. The control $u(t)$ [N] is a force applied to the first point mass $q_1(t)$. The vector $q(t)$ satisfies
\begin{equation}
\ddot{q}(t) + \mathbf{K} q(t) = \mathbf{f} u(t), \qquad \qquad q_i(0) = \tfrac{i-1}{10} - \tfrac{1}{2}, \qquad \dot{q}_i(0) = 0, \label{eq:dyn_q}
\end{equation}
where
\begin{equation}
\mathbf{K} = k\begin{bmatrix}
1 & -1 & 0 & 0 & \cdots & 0 & 0 & 0 \\
-1 & 2 & -1 & 0 & \cdots & 0 & 0 & 0 \\
0 & -1 & 2 & -1 & & 0 & 0 & 0 \\
0 & 0 & -1 & 2 & & 0 & 0 & 0 \\
\vdots & \vdots & & & \ddots & & & \vdots \\
0 & 0 & 0 & 0 & & 2 & -1 & 0 \\
0 & 0 & 0 & 0 & & -1 & 2 & -1 \\
0 & 0 & 0 & 0 & \cdots & 0 & -1 & 1 \\
\end{bmatrix}, \qquad 
\mathbf{f} = \begin{bmatrix}
1 \\ 0 \\ 0 \\ 0 \\ \vdots \\ 0 \\ 0 \\ 0
\end{bmatrix}. \label{eq:Kandf}
\end{equation}
Note that \eqref{eq:dyn_q} would be a (coarse) finite-difference discretization of the wave equation with Neumann boundary conditions if the  first and last rows in $\mathbf{K}$ are scaled by $2$, i.e.\ if the two masses at the end points $q_1(t)$ and $q_{11}(t)$ would have had mass $\tfrac{1}{2}$. This becomes a dynamical system of the form \eqref{eq:dyn_y} with $f(y,u) = Ay + Bu$ by setting 
\begin{equation}
y(t) = \begin{bmatrix}
q(t) \\ \dot{q}(t)
\end{bmatrix}, \qquad 
A = \begin{bmatrix}
O & I \\ -\mathbf{K} & O
\end{bmatrix}, \qquad 
B = \begin{bmatrix}
0 \\ \mathbf{f}
\end{bmatrix}. \label{eq:ex1_yAB}
\end{equation}
The RHC control for this system is computed by setting $E_T = 0$, $C = 10[I, O]$, and $R=1$ in \eqref{eq:JT}. The infinite-horizon optimal control $u_\infty^*(t)$ is obtained by solving the ARE \eqref{eq:ARE} using the MATLAB function \texttt{care} and discretizing \eqref{eq:xinf_dyn} by the Crank-Nicholson scheme. The optimal control problems on the finite horizon are solved by a gradient-descent algorithm using the Crank-Nicholson-based scheme from \cite{apel2012}. 

Figure \ref{fig:RHCnum} shows the obtained control $u_{\mathrm{RHC}}(t)$ for different values of $T$ and $\tau$. Figures \ref{fig:RHC1}, \ref{fig:RHC2}, and \ref{fig:RHC3} show that $u_{\mathrm{RHC}}(t)$ converges to $u^*_\infty(t)$ when $T - \tau$ is increased (for fixed $\tau = \tfrac{1}{2}$). Figure \ref{fig:RHC4} shows that decreasing $\tau$ while keeping $T - \tau$ fixed does not affect the RHC control $u_{\mathrm{RHC}}(t)$ visibly. These observations are in agreement with the estimates in Theorem \ref{thm:RHCconv} which only depend on $T - \tau$. Note that the dashed lines in Figure \ref{fig:RHCnum}, in Figure \ref{fig:RHC2} in particular, also show the parts of the controls $u^*_T(t,y_{\mathrm{RHC}}(k\tau), k\tau)$ that are not applied to the plant. Note that $u^*_T(k\tau + T,y_{\mathrm{RHC}}(k\tau), k\tau) = 0$ because the terminal cost $E_T = 0$. 

Figure \ref{fig:MPCAw} shows the influence of imperfections in the plant model on $u_{\mathrm{MPC}}(t)$. 

Figures \ref{fig:MPCAw1} and \ref{fig:MPCAw2} show the influence of a constant unit force applied to the rightmost mass $q_{11}(t)$, i.e.\ $f(y,u) = Ay + Bu$ and $w(t) = \bar{w} = [0_{1\times 21}, 1]^\top$. The MPC control $u_{\mathrm{MPC}}(t)$ in Figures \ref{fig:MPCAw1} and \ref{fig:MPCAw2} is compared to the optimal control for the infinite-horizon problem $u^*_\infty(t)$ (obtained from \eqref{eq:dyn_yinf} and \eqref{eq:dyn_xi}) and the limiting control for the MPC strategy $v_\infty(t)$ (obtained from \eqref{eq:dyn_zinf} and \eqref{eq:vinf}). Figures \ref{fig:MPCAw1} and \ref{fig:MPCAw2} indicate that decreasing $\tau$ brings $u_{\mathrm{MPC}}(t)$ closer to $v_\infty(t)$ (when $T - \tau$ is large enough), which is in agreement with Theorem \ref{thm:MPCconv}. 

Figures \ref{fig:MPCAw3} and \ref{fig:MPCAw4} show the control $u_{\mathrm{MPC}}(t)$ obtained when $w(t) \equiv 0$ and $f(y,u) = \tilde{A}y + Bu$ with
\begin{equation}
\tilde{A} = \begin{bmatrix}
O & I \\ - \mathbf{K} & 0.3I
\end{bmatrix}. \label{eq:tildeA_example}
\end{equation}
Note $\tilde{A}$ is not Hurwitz. 
The MPC control $u_{\mathrm{MPC}}(t)$ in Figures \ref{fig:MPCAw3} and \ref{fig:MPCAw4} is again compared to the infinite-horizon optimal control $u_\infty^*(t)$ (obtained based on the solution of the ARE \eqref{eq:ARE} with $A$ replaced by $\tilde{A}$), and the limiting MPC control $v_\infty(t)$ (obtained from \eqref{eq:dyn_zinf} and \eqref{eq:vinf}). Figures \ref{fig:MPCAw3} and \ref{fig:MPCAw4} indicate that reducing $\tau$ brings $u_{\mathrm{MPC}}(t)$ closer to $v_\infty(t)$ (when $T - \tau$ is sufficiently large), as Theorem \ref{thm:MPCconv} indicates. 

%The convergence rates predicted by Theorems \ref{thm:RHCconv} and \ref{thm:MPCconv} have been validated in a second numerical example in Section \ref{sec:example2} of the supplementary materials \cite{veldmanMPCsuppl2022}. 

\subsection{Example 2}
The convergence rates predicted by Theorems \ref{thm:RHCconv} and \ref{thm:MPCconv} are validated in a second example. The main motivation for considering a second example is that $\mu_\infty = 0.015\ldots$ in the example from Section \ref{sec:numerics}. This means that a large prediction horizon $T$ is required to make $e^{-\mu_\infty (T-\tau)}$ small, which makes the validation process computationally demanding. 

Therefore, a second numerical example is considered in which
\begin{equation}
A = -\mathbf{K}, \qquad 
B = \mathbf{f}, \qquad
E_T = 0, \qquad
C = I, \qquad
R = 1,
\end{equation}
where $\mathbf{K}$ and $\mathbf{f}$ are as in \eqref{eq:Kandf} with $k = 100$. The initial condition is chosen as $y_i(0) = \tfrac{i-1}{10} - \tfrac{1}{2}$.  Note that this example would be a (coarse) finite-difference discretization of the heat equation with Neumann boundary conditions if the  first and last rows in $\mathbf{K}$ are scaled by $2$. In this example, $\mu_\infty = 0.30\ldots$.

\begin{figure}
\begin{subfigure}{\columnwidth}
\centering
\includegraphics[width=0.55\textwidth]{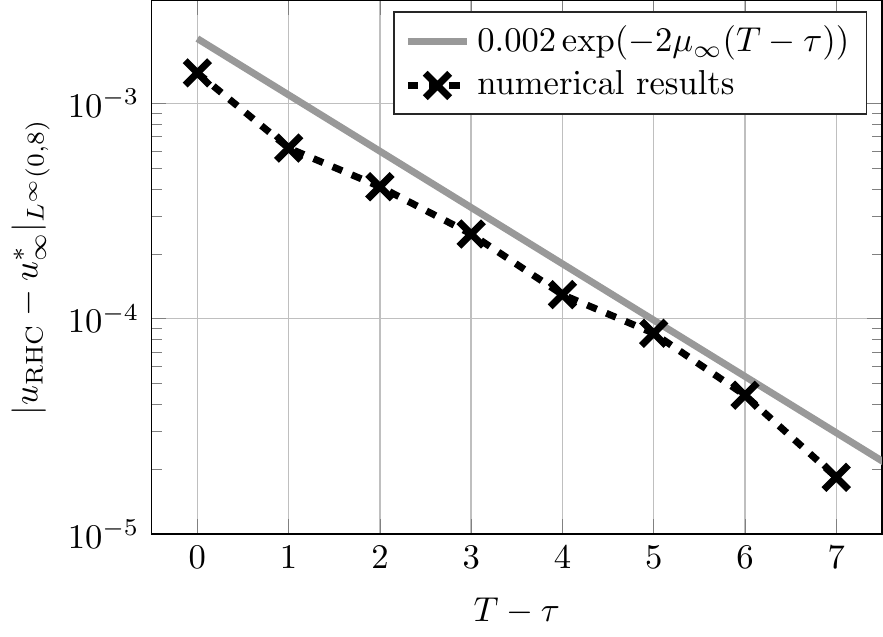}
\caption{The difference between $u_{\mathrm{RHC}}(t)$ and $u^*_\infty(t)$ for different values of $T-\tau$}
\label{fig:numval1}
\end{subfigure}
\begin{subfigure}{\columnwidth}
\centering
\includegraphics[width=0.55\textwidth]{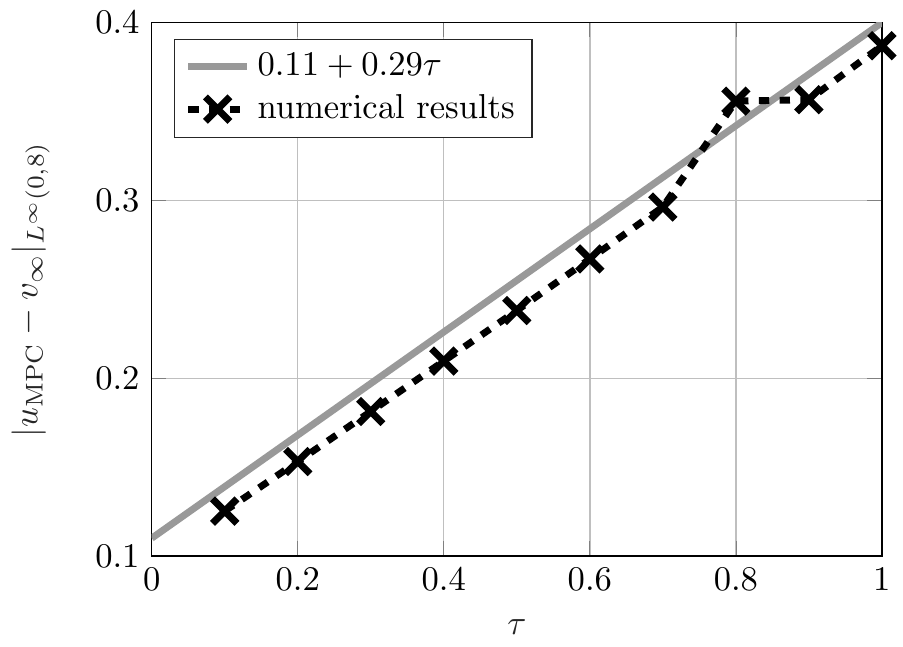}
\caption{The difference between $u_{\mathrm{MPC}}(t)$ and $v_\infty(t)$ with an additive disturbance $w(t) = \bar{w}$ for different values of $\tau$ with $T - \tau = 4$ fixed.}
\label{fig:numval2}
\end{subfigure}
\begin{subfigure}{\columnwidth}
\centering
\includegraphics[width=0.55\textwidth]{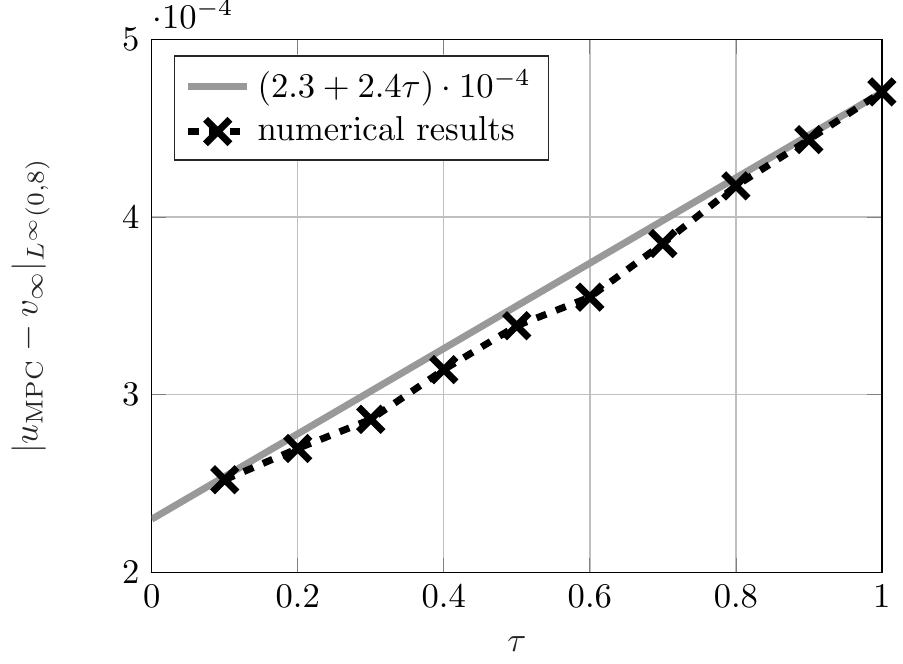}
\caption{The difference between $u_{\mathrm{MPC}}(t)$ and $v_\infty(t)$ when $f(y,u) = \tilde{A}y + Bu$, where $\tilde{A} = -\mathbf{K} + 0.3I$,  for different values of $\tau$ with $T - \tau = 4$ fixed.}
\label{fig:numval3}
\end{subfigure}
\caption{Numerical validation of the convergence rates for RHC and MPC derived in this paper for the second numerical example. }
\label{fig:numval}
\end{figure}

Figure \ref{fig:numval} shows that the convergence rates predicted by Theorems \ref{thm:RHCconv} and \ref{thm:MPCconv} can also be observed clearly in this numerical example. In particular, Figure \ref{fig:numval1} validates the convergence rate for RHC from Theorem \ref{thm:RHCconv} and shows that the difference between $u_{\mathrm{RHC}}$ and $u^*_\infty$ is proportional to $e^{-\mu_\infty(T-\tau)}$. Furthermore, Figures \ref{fig:numval2} and \ref{fig:numval3} validate the convergence rate for MPC from Theorem \ref{thm:MPCconv} and show that the difference between $u_{\mathrm{MPC}}$ and $v_\infty$ (obtained from \eqref{eq:dyn_zinf} and \eqref{eq:vinf}) is an affine function of $\tau$ (for $T-\tau$ fixed) when there are modeling errors $w(t)$ or $f(y,u) - Ay - Bu$. In particular, Figure \ref{fig:numval2} shows that this is the case for a constant disturbance $w(t) = \bar{w} = [0_{1\times 10}, 1]^\top$ and Figure \ref{fig:numval3} shows that this is the case when $f(y,u) = \tilde{A}y + Bu$, where $\tilde{A} = A + 0.3I$.

\section{Conclusions and discussions} \label{sec:conc}
%\subsection{Conclusions and discussions}

The results in this paper demonstrate that the local stability and convergence of RHC and MPC based on a linear plant model can be understood from the convergence of solution of the Riccati Differential Equation (RDE) to the (symmetric positive-definite) solution of the Algebraic Riccati Equation (ARE). The obtained estimates clearly show the influence of the two critical parameters in a MPC algorithm: the prediction horizon $T$ and the control horizon $\tau$. In particular, the optimal state trajectories and controls generated by RHC (i.e.\ MPC based on a perfect plant model) are close to their counterparts in an infinite-horizon optimal control problem when $T-\tau$ is sufficiently large. When the difference between the linear plant model and the nonlinear plant is sufficiently small, choosing the control horizon $\tau$ smaller reduces the influence of imperfections in the plant model. The obtained error estimates have been validated by numerical experiments.  

Several points deserve further discussion:
\begin{itemize} 
\item[1.] \textbf{Discrete-time systems} The results have been presented in a continuous-time setting, but have analogues in the discrete-time setting that is commonly used in the literature. Discrete-time versions of Lemma \ref{lem:Pconv} and the stability and convergence results for RHC are presented in Appendix \ref{sec:discretetime}. 

\item[2.] \textbf{Influence of the terminal cost}
The method in this paper does not require terminal constraints or terminal costs in the finite horizon optimal control problems and thus differs from several (mainly older) approaches in the existing literature. However, the estimates indicate that imposing a terminal cost $0 \preccurlyeq E_T \prec P_\infty$ improves the stability and convergence of the MPC strategy, see, e.g., \cite{allgower1999,ito2002} and Remark \ref{rem:K0ET} on page \pageref{rem:K0ET}.

\item[3.] \textbf{Estimates in existing literature}
The stability and convergence estimates in this paper depend in a different way on the control horizon $\tau$ than the estimates in the existing literature \cite{reble2012,azmi2016,azmi2018}, which cannot guarantee the stability of MPC for arbitrary small control horizons $\tau$, see Remark \ref{rem:existing} on page \pageref{rem:existing}. In contrast, the analysis in this paper shows that reducing the control horizon $\tau$ improves stability and convergence of the MPC strategy. Note however that the results in \cite{reble2012,azmi2016,azmi2018} also apply to MPC based on nonlinear plant models, while the approach in this paper is limited to MPC based on linear plant models. 

\item[4.] \textbf{Relation to the turnpike property} 
The analysis of MPC in this paper is based on the convergence of the solution of the RDE to the solution of the ARE. This is also the fundamental ingredient in one of the proofs of the turnpike in linear-quadratic optimal control, see \cite{porretta2013}. This indicates that there is an intimate relation between MPC and the turnpike property, which was also observed in \cite{grune2013, grune2020, pan2021}. 

\item[5.] \textbf{Computational advantage}  MPC only brings a computational advantage when the feedback operator for the infinite-horizon optimal control problem  cannot be computed easily, i.e.\ when the state space is high-dimensional or for constrained and/or nonlinear optimal control problems. 

\item[6.] \textbf{Infinite-dimensional systems} The extension of our results to infinite-dimensional problem should be relatively straightforward. In particular, infinite dimensional versions of Lemma \ref{lem:Pconv} have for example been obtained in \cite{porretta2013} in the context of turnpike and the Lipschitz condition for $f(y,u) - Ay - Bu$ in \eqref{eq:Lipschitz} could be relaxed to a monotonicity condition that is applicable to semilinear Partial Differential Equations (PDEs). 

\item[7.] \textbf{Constrained and/or nonlinear plant models}
The results in this paper have been derived under the limiting assumption that the plant model used in the MPC controller is Linear Time Invariant (LTI) and that there are no constraints on the control or state of the plant to be controlled. A natural way to extend the approach from this paper to constrained and/or nonlinear plant models is through Hamilton-Jacobi theory, which has also been applied in the context of turnpike, see, e.g., \cite{anderson1987}. The main approach and remaining problems are outlined in Appendix \ref{sec:HJ}.

\item[8.] \textbf{Adaptive MPC}
Theorem \ref{thm:MPCconv} indicates that decreasing $\tau$ improves the robustness of the MPC strategy against modeling errors. Because decreasing $\tau$ also increases the computational cost for the MPC strategy, it is natural to implement MPC on an adaptive timegrid $0 = t_0 < t_1 < t_2 < \ldots$, see e.g.\ \cite{grune2010}. For $T$  sufficiently large, Theorem \ref{thm:MPCconv} suggests that the control horizon $\tau_k = t_k - t_{k-1}$ should be small when $|w|_{L^\infty(t_{k-1},t_k)}$ is large and that $\tau_k$ can be increased when $|w|_{L^\infty(t_{k-1},t_k)}$ is small. %MPC with an adaptive prediction horizon $T = T_k$ as for example considered in \cite{sun2021} is of interest. 
Although MPC with adaptive prediction and/or control horizons has been proposed before (see, e.g., \cite{kim2010, sun2021}), the error estimates in this paper may lead to new insights. 

\item[9.] \textbf{Deep learning}
Because the training of Deep Neural Networks (DNNs) can be viewed as a (nonlinear) optimal control problem, see e.g.\ \cite{E2017, benning2019, esteve2021, esteve2021sparse}, the ideas of MPC can also be applied in this context. In particular, instead of training all layers in a deep neural network simultaneously, one can adopt a receding horizon approach and first train only the first $T$ layers, fix the found weights in the first $\tau \leq T$ layers, shift the considered time horizon by $\tau$ and repeat. The results in this paper are thus also of interest for deep learning. 

\end{itemize}

\appendix
\section{Long-term behavior of the RDE} \label{app:RDE}

For completeness, the proof of Lemma \ref{lem:Pconv} is given below. 
The proof is inspired by similar results in \cite{callier1994, porretta2013}. 

\begin{proof}
Write $\mathcal{E}(t) := P_\infty - \mathcal{P}(t)$ and substract \eqref{eq:RDE_reversed} from \eqref{eq:ARE} to find 
\begin{align}
\dot{\mathcal{E}}(t) &= A^\top \mathcal{E}(t) + \mathcal{E}(t)A - P_\infty B R^{-1} B^\top P_\infty + \mathcal{P}(t) B R^{-1} B^\top \mathcal{P}(t) \nonumber \\
&= A^\top \mathcal{E}(t) + \mathcal{E}(t)A - \mathcal{E}(t) B R^{-1} B^\top P_\infty - (P_\infty - \mathcal{E}(t)) B R^{-1} B^\top \mathcal{E}(t) \nonumber \nonumber \\
&= A_\infty^\top \mathcal{E}(t) + \mathcal{E}(t) A_\infty +  \mathcal{E}(t) B R^{-1} B^\top \mathcal{E}(t), \label{eq:dyn_E}
\end{align}
where the last equality follows from the definition of $A_\infty$ in \eqref{eq:Ainfty}. The IC in \eqref{eq:RDE_reversed} also shows that $\mathcal{E}(0) = P_\infty - E_T$. 

Next, introduce the function $\mathcal{S}(t)$ by the relation
\begin{equation}
\mathcal{E}(t) = e^{A_\infty^\top t} \mathcal{S}(t) e^{A_\infty t}. \label{eq:def_E}
\end{equation}
Inserting $\dot{\mathcal{E}}(t) = A_\infty^\top \mathcal{E}(t) + \mathcal{E}(t) A + e^{A_\infty^\top t} \dot{\mathcal{S}}(t) e^{A_\infty t}$ into \eqref{eq:dyn_E}, shows that
\begin{equation}
\dot{\mathcal{S}}(t) = \mathcal{S}(t) e^{A_\infty t} B R^{-1} B^\top e^{A_\infty^\top t} \mathcal{S}(t). \label{eq:dyn_S}
\end{equation}
Note that  $\dot{\mathcal{S}}(t) \succcurlyeq 0$, so that $\mathcal{S}(0)\preccurlyeq \mathcal{S}(t) \preccurlyeq \mathcal{S}(\infty)$ for all $t \geq 0$ which implies that 
\begin{equation}
\| \mathcal{S}(t) \| \leq \max\{ \| \mathcal{S}(0) \|, \| \mathcal{S}(\infty) \| \}. \label{eq:Sbound}
\end{equation}
Taking norms in \eqref{eq:def_E} using \eqref{eq:growthbound} and \eqref{eq:Sbound} now shows that
\begin{equation}
\| \mathcal{P}(t) - P_\infty\| \leq M_\infty^2 \max\{ \| \mathcal{S}(0) \|, \| \mathcal{S}(\infty) \| \} e^{-2 \mu_\infty t}. \label{eq:est_E}
\end{equation}
This is an estimate of the form \eqref{eq:Pconv}. 
Because $\mathcal{S}(0) = \mathcal{E}(0) = P_\infty - E_T$, it remains to compute $\mathcal{S}(\infty)$. To this end, introduce
\begin{align}
\mathcal{W}(t) &:= \int_0^t e^{A_\infty s} B R^{-1} B^\top e^{A_\infty^\top s} \ \mathrm{d}s.   \label{eq:def_W}\\ 
\tilde{\mathcal{S}}(t) &:= \mathcal{S}(0) (I - \mathcal{W}(t)\mathcal{S}(0))^{-1}. \label{eq:def_tildeS}
\end{align}
We claim that $\tilde{\mathcal{S}}(t) = \mathcal{S}(t)$. By definition, $\tilde{S}(0) = \mathcal{S}(0)$. Furthermore, if $I - \mathcal{W}(t)\mathcal{S}(0)$ is invertible, differentiating \eqref{eq:def_tildeS} shows that $\dot{\tilde{\mathcal{S}}}(t) = \tilde{\mathcal{S}}(t) \dot{\mathcal{W}}(t) \tilde{\mathcal{S}}(t)$. Computing $\dot{\mathcal{W}}(t)$ using \eqref{eq:def_W} shows that $\tilde{\mathcal{S}}(t)$ satisfies \eqref{eq:dyn_S}. So if $I - \mathcal{W}(t)\mathcal{S}(0)$ is invertible for all $t$, 
\begin{equation}
\mathcal{S}(\infty) = \tilde{\mathcal{S}}(\infty) = \mathcal{S}(0) (I - \mathcal{W}(\infty)\mathcal{S}(0))^{-1}.
\end{equation}
Comparing \eqref{eq:est_E} and \eqref{eq:Pconv} and using that $\mathcal{S}(0) = P_\infty - E_T$, it follows that the constant $K_0$ in \eqref{eq:Pconv} can be chosen as
\begin{equation}
K_0 = M_\infty^2 \max\{\| P_\infty - E_T \|, \| (P_\infty - E_T)(I - \mathcal{W}(\infty)(P_\infty - E_T))^{-1} \| \}. \label{eq:def_K0}
\end{equation}

It remains to show that $I - \mathcal{W}(t)\mathcal{S}(0)$ is invertible for all $t$. To this end, note that $\mathcal{W}(t)$ is a (weighted closed-loop) controllability Grammian and is thus invertible for all $t$ because $(A,B)$ is controllable. Therefore,
\begin{equation}
I - \mathcal{W}(t)\mathcal{S}(0) = \mathcal{W}(t)((\mathcal{W}(t))^{-1} - P_\infty + E_T).
\end{equation}
Because $E_T \succcurlyeq 0$ and $\dot{\mathcal{W}}(t) \succcurlyeq 0$
\begin{equation}
(\mathcal{W}(t))^{-1} - P_\infty + E_T \succcurlyeq (\mathcal{W}(\infty))^{-1} - P_\infty. 
\end{equation}
It thus suffices to show that $(\mathcal{W}(\infty))^{-1} - P_\infty$ is positive definite. To see this, note that the definition of $\mathcal{W}(t)$ in \eqref{eq:def_W} implies that $\mathcal{W}(\infty)$ is the solution of the Lyapunov equation
\begin{equation}
A_\infty \mathcal{W}(\infty) + \mathcal{W}(\infty) A_\infty^\top + BR^{-1}B^\top = 0. 
\end{equation}
Multiplying this equation from both sides by $(\mathcal{W}(\infty))^{-1}$ and subtracting the result from the ARE \eqref{eq:ARE}, it can be shown that $Q_\infty := P_\infty - (\mathcal{W}(\infty))^{-1}$ is also a solution of the ARE \eqref{eq:ARE}. As having $Q_{\infty} = P_{\infty}$ implies that $(\mathcal{W}(\infty))^{-1} = 0$ (which is absurd), it follows that $Q_{\infty}$ is the unique symmetric \emph{negative-definite} solution of \eqref{eq:ARE}. Therefore, $(\mathcal{W}(\infty))^{-1}  - P_\infty = -Q_\infty$ is positive definite and $I - \mathcal{W}(t) \mathcal{S}(0)$ is invertible for all time. 
\end{proof}

\section{Hamilton-Jacobi theory} \label{sec:HJ}
This paper has focused on the analysis of MPC based on a linear unconstrained plant model. Hamilton-Jacobi theory shows some potential to overcome this limitation. The main approach and the remaining problems are discussed in this appendix. 

Consider a RHC strategy for the plant \eqref{eq:dyn_y} with $w(t) \equiv 0$ and the additional requirement that the control $u(t)$ should take values in a nonempty closed and convex set $U_{\mathrm{ad}} \subseteq \mathbb{R}^m$. Algorithm \ref{alg:MPC} remains essentially unchanged, but the control $u^*_T(t;x_1,t_1)$ is now computed as the minimizer of $J_T(u;x_1, t_1)$
\begin{equation}
I_T(u) = \frac{1}{2}(x(t_1+T))^\top E_T x(t_1+T)) \\
+ \frac{1}{2}\int_{t_1}^{t_1+T} | Cx(t)|^2 + (u(t))^\top Ru(t) \ \mathrm{d}t, \label{eq:JT_general}
\end{equation}
over all $u \in L^2(0,T; U_{\mathrm{ad}})$ subject to the dynamics
\begin{equation}
\dot{x}(t) = f(x(t), u(t)), \qquad \qquad x(t_1) = x_1, \label{eq:x_general}
\end{equation}
with $f$ as in \eqref{eq:dyn_y}, $C \in \mathbb{R}^{\ell \times n}$ and $R \succ 0$ as in $I_\infty(u)$ in \eqref{eq:def_I}, and a terminal cost $E_T \succcurlyeq 0$. Key quantities in the following discussion are the value functions
\begin{align}
V_\infty(x_1) = \inf_{L^2(0,\infty;U_{\mathrm{ad}})} I_\infty(u;x_1,t_1), \qquad
 V_T(x_1) = \inf_{L^2(0,T;U_{\mathrm{ad}})} I_T(u;x_1,t_1). \label{eq:def_V}
\end{align}
Note that $V_\infty(x_1)$ and $V_T(x_1)$ do not depend on $t_1$ because the problem \eqref{eq:JT_general}--\eqref{eq:x_general} is time invariant. For simplicity, it is assumed that $V_\infty(x_1)$ is finite for all $x_1$, which is true when \eqref{eq:x_general} is null controllable.

It is well-known that $V_\infty(x_1)$ and $V_T(x_1)$ are differentiable almost everywhere, see, e.g., \cite{lions1982}. If $V_\infty(x_1)$ and $V_T(x_1)$ are differentiable everywhere, the optimal state trajectories $x^*_\infty(t)$ and $x^*_T(t)$ are given by a feedback law, i.e.\ they satisfy
\begin{equation}
\dot{x}^*_\infty(t) = f(x^*_\infty(t)), U_\infty(x^*_\infty(t))), \qquad \dot{x}^*_T(t) = f(x^*_\infty(t)), U_T(t,x^*_\infty(t))), 
\end{equation}
where
\begin{align}
U_\infty(x) &= \mathrm{arg} \min_{u\in U_{\mathrm{ad}}} \left\{ \frac{\mathrm{d}V_\infty}{\mathrm{d}x_1}(x) f(x,u) + |Cx|^2 + u^\top R u \right\}, \label{eq:Uinf_feedback} \\
U_T(t,x) &= \mathrm{arg} \min_{u\in U_{\mathrm{ad}}} \left\{ \frac{\partial V_{T-t}}{\partial x_1}(x) f(x,u) + |Cx|^2 + u^\top R u \right\}. \label{eq:UT_feedback}
\end{align}

It is now easy to see that $y_{\mathrm{RHC}}(t) \rightarrow x^*_\infty(t)$ if $U_T(t,x) \rightarrow U_\infty(x)$ for $T \rightarrow \infty$ (and $t \leq \tau$ fixed). According to \eqref{eq:Uinf_feedback} and \eqref{eq:UT_feedback}, the latter condition is satisfied when
\begin{equation}
\lim_{T\rightarrow \infty} \frac{\partial V_T}{\partial x_1}(x_1) = \frac{\mathrm{d}V_\infty}{\mathrm{d}x_1}(x_1). \label{eq:dVdx_conv}
\end{equation}
When an explicit error bound for \eqref{eq:dVdx_conv} can be obtained, explicit stability and convergence conditions for RHC and MPC based on a nonlinear and constrained plant model can be obtained along the lines of this paper. However, because $V_T$ and $V_\infty$ are generally not differentiable everywhere, verifying \eqref{eq:dVdx_conv} is not trivial. 

\begin{remark} \label{rem:LTIcase}
For the linear unconstrained system model considered in this paper, i.e.\  for $f(y) = Ay$ and $U_{\mathrm{ad}} = \mathbb{R}^m$, $V_\infty(x_1) = \tfrac{1}{2}x_1^\top P_\infty x_1$ and $V_T(x_1) = \tfrac{1}{2}x_1^\top \mathcal{P}(T) x_1$ where $P_\infty$ is the symmetric positive definite solution of the ARE \eqref{eq:ARE} and $\mathcal{P}(t)$ is the solution of time-reversed RDE \eqref{eq:RDE_reversed}, see e.g.\ \cite{sage1968}. Lemma \ref{lem:Pconv} is thus assures that \eqref{eq:dVdx_conv} holds in the unconstrained linear quadratic case. 
\end{remark}

\section{Discrete time RHC} \label{sec:discretetime}
Because a large part of the literature on MPC focusses on discrete-time systems, the discrete-time analogues of the results for RHC from Subsection \ref{ssec:RHC} are proved in this appendix. 

\subsection{Discrete-time RHC}
In discrete-time RHC, the goal is to control the linear dynamics
\begin{equation}
x_{t+1} = Ax_t + Bu_t. \label{eq:dyn_x_DT}
\end{equation}
Just as in the continuous time setting, $x_t$ evolves in $\mathbb{R}^n$, $x_0 \in \mathbb{R}^n$ is a given initial condition,  $A \in \mathbb{R}^{n \times n}$ is the system matrix, and $B \in \mathbb{R}^{n \times m}$ is the input matrix, but the time $t$ now takes discrete values $t \in \mathbb{N}$. It is again assumed that the state $x_t$ can be measured at certain time instances $t = k\tau$ for some $\tau \in \mathbb{N}$ fixed and $k \in \mathbb{N}$. 

Similarly as before, introduce 
\begin{equation}
u^T(\bar{x},t_1) = (u^T_{t_1}(\bar{x},t_1), u^T_{t_1+1}(\bar{x},t_1), \ldots, u^T_{t_1+T-1}(\bar{x},t_1)),
\end{equation} 
as the minimizer of
\begin{equation}
J^T(u,\bar{x},t_1) = \frac{1}{2}x_{t_1+T}^\top E_T x_{t_1+T} + \frac{1}{2} \sum_{t=t_1}^{t_1+T-1} \left( |Cx_t|^2  +  u_t^\top R u_t \right), \label{eq:JT_DT}
\end{equation}
subject to the dynamics
\begin{equation}
x_{t+1} = Ax_t + Bu_t, \qquad \qquad x_{t_1} = \bar{x}. \label{eq:dyn_xT_DT}
\end{equation}
Here, $E_T$ is a symmetric positive semi-definite matrix, $R \in \mathbb{R}^{m\times m}$ is a symmetric and positive definite matrix, $C \in \mathbb{R}^{\ell \times  n}$ (with $\ell \leq n$) is the output matrix. It is well-known that $u^T(\bar{x},t_1)$ exists an is unique (for  all $\bar{x}$ and $t_1$). 

The discrete-time RHC algorithm given in Algorithm \ref{alg:RHC_DT} is similar to Algorithm \ref{alg:MPC} for continuous time MPC. Note that the RHC control is denoted by $u^{\mathrm{RHC}}_t$ and the corresponding state trajectory by $x^{\mathrm{RHC}}_t$. 

\begin{algorithm}
\textbf{Step 1} Choose a prediction horizon $T \in \mathbb{N}$ and a control horizon $\tau \in \{1,2, \ldots, T \}$. Set $k = 0$. 

\textbf{Step 2} Measure $x^{\mathrm{RHC}}_{k\tau}$ and compute the control $u^T(x^{\mathrm{RHC}}_{k\tau}, k\tau)$.

\textbf{Step 3} Set $u^{\mathrm{RHC}}_t = u^T_t(x^{\mathrm{RHC}}_{k\tau}, k\tau)$ for $t \in \{k\tau, k\tau+1, \ldots, k\tau + \tau-1 \}$ and \\ $x^{\mathrm{RHC}}_t = x^T_t(x^{\mathrm{RHC}}_{k\tau}, k\tau)$ for $t \in \{k\tau+1, k\tau+2, \ldots, k\tau+\tau \}$. 

\textbf{Step 4} Increase $k$ by 1 and go to step 2. 

\caption{Receding Horizon Control for Discrete Time Systems}
\label{alg:RHC_DT}
\end{algorithm}
Again, the question arises for which prediction horizons $T$ and control horizons $\tau$ this control strategy is stabilizing. Just as in the continuous time setting, this question is closely related to the minimization of the infinite horizon cost
\begin{equation}
J^\infty(u; x_0) = \frac{1}{2} \sum_{t=0}^{\inf} \left( |Cx_t|^2  +  u_t^\top R u_t \right), \label{eq:Jinf_DT}
\end{equation}
subject to the dynamics
\begin{equation}
x_{t+1} = Ax_t + Bu_t, \qquad \qquad x_0 = x_0. \label{eq:dyn_xinf_DT}
\end{equation}
Note that the minimizer $u^\infty_t$ exists if $(A,B)$ is controllable. 
The relation between this optimal control problem and the RHC algorithm will be made more precise with the Riccati theory in the next subsection. 

\begin{remark}[Relation to time-discretization]
Applying a first-discretize-then-optimize approach (FDTO) with a fixed step size $\Delta t$ to the continuous-time RHC algorithm leads to a problem setting as described in this subsection. It is therefore clear that analogues of the continuous-time results will also hold in the discrete-time setting when $\Delta t$ is sufficiently small. 

The results in this section are stronger than this. They show that the discrete-time MPC controller stabilizes the discrete time-system for $T-\tau$ sufficiently large provided that the discrete-time system is controllable and observable, a condition that does not require $\Delta t$ to be small. 

Another natural question is how small $\Delta t$ needs to be such that the discrete-time MPC-controller stabilizes the continuous-time system and how this required $\Delta t$ depends on $T$ and $\tau$. This question can be answered using the robustness results from Theorems \ref{thm:stabRHC_DT} and \ref{thm:convRHC_DT} by viewing the time-discretization error in the control as an additive disturbance $w(t)$. Error bounds for the time-discretization of optimal control problems are available, see, e.g., \cite{hager2000, gerdts2012} and the references therein. 
\end{remark}

\subsection{Riccati theory}
It is well-known that the optimal state trajectory for the optimal control problem \eqref{eq:JT_DT}--\eqref{eq:dyn_xT_DT}
 can be computed as the solution of
\begin{equation}
x^T_{t+1}(\bar{x},t_1) = \left( I + B R^{-1} B^\top Q^T_{t+1-t_1} \right)^{-1} A x^T_t(\bar{x},t_1), \qquad  x^T_{t_1}(\bar{x},t_1) = \bar{x}.  \label{eq:xT_dyn_DT}
\end{equation}
Here, $Q^T_t$ is the solution of the Discrete-time Riccati Difference Equation (DRDE)
\begin{equation}
Q^T_t = A^\top Q^T_{t+1} (I + BR^{-1}B^\top Q^T_{t+1})^{-1} A + C^\top C, \qquad\qquad\quad Q^T_T = E_T, \label{eq:DRDE}
\end{equation}
which is again solved backward in time starting from the final condition. It is therefore convenient to introduce $\mathcal{Q}_t$ (for all $t \in \mathbb{N}$) as the solution of the time-reversed DRDE
\begin{equation}
\mathcal{Q}_{t+1} = A^\top \mathcal{Q}_t (I + B R^{-1} B^\top \mathcal{Q}_t )^{-1} A + C^\top C, \qquad\qquad \mathcal{Q}_0 = E_T. \label{eq:DRDEreversed}
\end{equation} 
Note that
\begin{equation}
Q^T_{t} = \mathcal{Q}_{T-t}. \label{eq:defQ_DT}
\end{equation}

Similarly, if $(A,B)$ is controllable and $(A,C)$ is observable, the optimal state trajectory for the infinite-horizon optimal control problem \eqref{eq:Jinf_DT}--\eqref{eq:dyn_xinf_DT} is 
\begin{equation}
x^\infty_{t+1} = \left( I + B R^{-1} B^\top \mathcal{Q}_\infty \right)^{-1} A x^\infty_t, \label{eq:x_inf_opt_DT}
\end{equation}
where $Q_\infty$ is the unique positive definite solution of the Discrete-time Algebraic Riccati Equation (DARE)
\begin{equation}
Q_\infty = A^\top Q_\infty (I + BR^{-1}B^\top Q_\infty)^{-1} A + C^\top C. \label{eq:DARE}
\end{equation}
To simplify notation, introduce
\begin{equation}
M_\infty := (I + BR^{-1}B^\top Q_\infty)^{-1}, \qquad M_t := (I + BR^{-1}B^\top \mathcal{Q}_t)^{-1}, \label{eq:def_MinfMt}
\end{equation}
Then \eqref{eq:xT_dyn_DT} and \eqref{eq:x_inf_opt_DT} can be rewritten as
\begin{equation}
x^T_{t+1} = M_{T-t}Ax^T_t, \qquad \qquad  x^\infty_{t+1} = M_\infty A x^\infty_t.
\end{equation}
The proof of convergence of the RHC strategy again relies on the convergence of $\mathcal{Q}_t$ to $Q_\infty$ for $t \rightarrow \infty$. An explicit error estimate is proved in the next subsection. 

\begin{remark}
Writing the Euler-Lagrange equations for the discrete-time optimal control problem in \eqref{eq:JT_DT} and \eqref{eq:dyn_xT_DT}, it follows that the optimal state trajectory satisfies (see, e.g., \cite[Example 6.2-1]{sage1968})
\begin{equation}
x_{t+1}^T = (I - BR^{-1}B^\top Q^T_{t+1}(I+BR^{-1}B^\top Q^T_{t+1})^{-1})A x_t^T, \qquad x_0^T = \bar{x}, \label{eq:xT_dyn_DT_1}
\end{equation}
where $Q^T_t$ is the symmetric positive-definite solution of the DRDE \eqref{eq:DRDE}. 
Because $I - X(I+X)^{-1} = (I + X)^{-1}$, \eqref{eq:xT_dyn_DT_1}  can be rewritten in the form \eqref{eq:xT_dyn_DT}. 
\end{remark}

\begin{remark}
Equations \eqref{eq:xT_dyn_DT} and \eqref{eq:DRDE} can be rewritten into their more commonly found form, using the Woodbury matrix identity, see e.g.\ \cite{kucera1972, henderson1981}, which states that for any invertible matrices $X \in \mathbb{R}^{n_1 \times n_1}$ and $Y \in \mathbb{R}^{n_2 \times n_2}$ and matrices $\tilde{B} \in \mathbb{R}^{n_1 \times n_2}$ and $\tilde{C} \in \mathbb{R}^{n_2 \times n_1}$
\begin{equation}
(X^{-1}+\tilde{B}Y^{-1}\tilde{C})^{-1} = X - X\tilde{B}(Y + \tilde{C}X\tilde{B})^{-1}\tilde{C}X, 
\end{equation}
it follows that (setting $X = I$, $Y=R$, $\tilde{B}=B$, and $\tilde{C} = B^\top Q^T_{t+1}$)
\begin{equation}
(I + BR^{-1}B^\top Q^T_{t+1})^{-1} = I - B(R + B^\top Q^T_{t+1} B)^{-1} B^\top Q^T_{t+1}.
\end{equation}
Inserting this into \eqref{eq:xT_dyn_DT} and \eqref{eq:DRDE} we obtain the more commonly found forms 
\begin{align}
x_{t+1}^T &= (I - B(R + B^\top Q^T_{t+1} B)^{-1}B^\top Q^T_{t+1}) Ax_t^T, \qquad\qquad x_0^T = \bar{x}, \label{eq:xT_dyn_DT_2} \\
Q^T_t &= A^\top Q^T_{t+1} A - A^\top Q^T_{t+1} B (R + B^\top Q^T_{t+1} B)^{-1} B^\top Q^T_{t+1} A + C^\top C, \label{eq:DRDE_1}
\end{align}
with the final condition $Q^T_T = E_T$. 
\end{remark}

\subsection{Long-term behavior of the DRDE}
The main result of this subsection is the following. 

\begin{lemma} \label{lem:Qconv}
Let $Q_\infty$ be the symmetric positive definite solution of \eqref{eq:DARE}, $\mathcal{Q}_t$ be the solution of \eqref{eq:DRDEreversed} with $0 \preccurlyeq E_T \prec Q_\infty$, then there exists a constant $\tilde{K}_0$ such that
\begin{align}
\| \mathcal{Q}_t - Q_\infty \| &\leq \tilde{K}_0 \| M_\infty A \|^{2t}. \label{eq:lem_Qconv}
\end{align}
\end{lemma}
Note that $M_\infty A$ is the matrix that generates the closed-loop dynamics \eqref{eq:x_inf_opt_DT}. Because the minimal infinite horizon cost is finite and $(A,C)$ is observable, it follows that $\| M_\infty A \| < 1$. The proof uses the following basic result from linear algebra. 
\begin{lemma} \label{lem:linalg}
Let $X,Y \in \mathbb{R}^{n \times n}$ be symmetric positive semi-definite matrices, then 
\begin{itemize}
\item[(i)] $I + XY$ and $I + YX$ are invertible,
\item[(ii)] $(I+XY)^{-1}X$ and $Y(I+XY)^{-1}$ are symmetric, and
\item[(iii)] $0 \preccurlyeq (I + XY)^{-1}X \preccurlyeq X$ and $0 \preccurlyeq Y (I+XY)^{-1} \preccurlyeq Y$. 
\end{itemize}
\end{lemma}
\begin{proof} As $X$ is symmetric positive semidefinite, there exists a symmetric positive semi definite matrix $X^{1/2}$ such that $X = X^{1/2}X^{1/2}$. Because the eigenvalues of $\tilde{X}\tilde{Y}$ are the same as the eigenvalues of $\tilde{Y}\tilde{X}$ for all matrices $\tilde{X}$ and $\tilde{Y}$, the eigenvalues of $X^{1/2} Y X^{1/2}$, $X^{1/2}X^{1/2} Y = XY$, and $Y X^{1/2} X^{1/2} = YX$ are the same. Since and $X^{1/2}YX^{1/2}$ is clearly positive semi-definite, all eigenvalues of $XY$ and $YX$ are nonnegative. Therefore, all eigenvalues of $I + XY$ and $I + YX$ are bigger than 1 and the matrices are invertible. 

For point (ii), note that (for any matrices $X$ and $Y$ for which $I+XY$ and $I+YX$ are invertible)
\begin{equation}
(I+XY)^{-1}X = X (I+YX)^{-1}, %= \left( (I + X^\top Y^\top)^{-1} X^\top \right)^\top = \left( (I + X Y)^{-1} X \right)^\top,
\label{eq:lem_linalg_step1}
\end{equation}
which can be verified by  multiplying by $(I + XY)$ from the left and by $(I + YX)$ from the right, see also, e.g.,  \cite{henderson1981}. Point (ii) now follows by computing the transpose of the expression on the right, using that $X$ and $Y$ are symmetric.

For point (iii), note that replacing $(X,Y)$ in \eqref{eq:lem_linalg_step1} by $(X^{1/2}, X^{1/2}Y)$ yields
\begin{equation}
(I + XY)^{-1}X^{1/2} = X^{1/2}(I + X^{1/2}YX^{1/2})^{-1}.
\end{equation}
Since $(I + X^{1/2}YX^{1/2})^{-1} \preccurlyeq I$, multiplying this equation from the left by $x^\top$ and from the right by $X^{1/2}x$ shows that
\begin{equation}
0 \leq x^\top (I+XY)^{-1}X x = x X^{1/2}(I + X^{1/2}YX^{1/2})^{-1} X^{1/2}x \leq x^\top X x, 
\end{equation}
for all vectors $x \in \mathbb{R}^n$ and the result follows. 
\end{proof}

It is now possible to prove Lemma \ref{lem:Qconv}. 

\begin{proof} 
With the definitions of $M_\infty$ and $M_t$ in \eqref{eq:def_MinfMt}, the DARE \eqref{eq:DARE} and and DRDE \eqref{eq:DRDEreversed} can be rewritten as
\begin{equation}
Q_\infty = A^\top Q_\infty M_\infty A, \qquad \qquad \mathcal{Q}_{t+1} = A^\top \mathcal{Q}_t M_t A. \label{eq:lem_Qconv_step2}
\end{equation}
Writing $\mathcal{E}_t := Q_\infty - \mathcal{Q}_t$, it follows that $\mathcal{E}_0 = Q_\infty - E_T$ and that
\begin{align}
\mathcal{E}_{t+1} &= A^\top Q_\infty M_\infty A - A^\top \mathcal{Q}_t M_t A \\
&= A^\top \mathcal{E}_t M_\infty A + A^\top \mathcal{Q}_t (M_\infty - M_t)A \nonumber \\
&= A^\top \mathcal{E}_t M_\infty A + A^\top \mathcal{Q}_t M_t (M_t^{-1} - M_\infty^{-1}) M_\infty A \nonumber \\
&= A^\top \mathcal{E}_t M_\infty A - A^\top \mathcal{Q}_t M_t B R^{-1} B^\top \mathcal{E}_t M_\infty A \nonumber \\
&= A^\top (I - \mathcal{Q}_t M_t B R^{-1} B^\top ) \mathcal{E}_t M_\infty A \nonumber \\
&= A^\top \underbrace{M_\infty^\top (I + Q_\infty B R^{-1} B^\top)}_{=I} (I - \mathcal{Q}_t M_t B R^{-1} B^\top ) \mathcal{E}_t M_\infty A, \label{eq:lem_Qconv_step3}
\end{align}
where the fourth and the last equality follow from \eqref{eq:def_MinfMt}. 
Expanding the brackets for the factor in the middle of the last expression yields
\begin{align}
(I + Q_\infty & B R^{-1} B^\top) (I - \mathcal{Q}_t M_t B R^{-1} B^\top ) \nonumber \\
&= I + Q_\infty B R^{-1}B^\top - \mathcal{Q}_t M_t B R^{-1} B^\top - Q_\infty B R^{-1} B^\top \mathcal{Q}_t M_t B R^{-1} B^\top \nonumber \\
&= I + Q_\infty \underbrace{(I + BR^{-1}B^\top \mathcal{Q}_t )M_t}_{=I} BR^{-1}B^\top \nonumber \\
&\qquad\qquad - \mathcal{Q}_t M_t B R^{-1} B^\top - Q_\infty B R^{-1} B^\top \mathcal{Q}_t M_t B R^{-1} B^\top \nonumber \\
&= I + \mathcal{E}_t M_t B R^{-1} B^\top, \label{eq:lem_Qconv_step4}
\end{align}
so that \eqref{eq:lem_Qconv_step3} can be written as
\begin{align}
\mathcal{E}_{t+1} 
&= A^\top M_\infty^\top (I + \mathcal{E}_tM_tBR^{-1}B^\top) \mathcal{E}_t M_\infty A \nonumber \\ 
&= A^\top M_\infty^\top (\mathcal{E}_t + \mathcal{E}_tM_tBR^{-1}B^\top \mathcal{E}_t) M_\infty A  \label{eq:lem_Qconv_step5}
\end{align}
Now introduce a new variable $\mathcal{S}_t$ by the relation
\begin{equation}
\mathcal{E}_t = (A^\top M_\infty^\top)^t \mathcal{S}_t (M_\infty A)^t. \label{eq:lem_Qconv_step6}
\end{equation}
Inserting this expression into \eqref{eq:lem_Qconv_step5} shows that $\mathcal{S}_t$ should satisfy
\begin{equation}
\mathcal{S}_{t+1} = \mathcal{S}_t + \mathcal{S}_t (M_\infty A)^t M_t B R^{-1} B^\top (A^\top M_\infty^\top)^t \mathcal{S}_t =: \mathcal{S}_t + \mathcal{S}_t Z_t \mathcal{S}_t,\label{eq:lem_Qconv_step7}
\end{equation}
where the matrix $Z_t$ has been introduced for brevity. Applying Lemma \ref{lem:linalg} with $X = BR^{-1}B^\top$ and $Y = \mathcal{Q}_t$ shows that $M_t B R^{-1}B^\top$ is symmetric positive semidefinite. Therefore $Z_t$ is also symmetric and positive semi definite and $\mathcal{S}_t$ is nondecreasing. Because $E_T \prec Q_\infty$ by assumption, $0\prec \mathcal{S}_0 \preccurlyeq \mathcal{S}_t \preccurlyeq \mathcal{S}_\infty$ and 
\begin{equation}
\| \mathcal{S}_t \| \leq \| \mathcal{S}_\infty \|. 
\end{equation}

It thus remains to find an upper bound for $\mathcal{S}_\infty$, which is equivalent to a lower bound for $\mathcal{S}_\infty^{-1}$.  Now observe that
\begin{align}
\mathcal{S}_{t+1}^{-1} &= (\mathcal{S}_t(I+Z_t\mathcal{S}_t))^{-1} = (I + Z_t \mathcal{S}_t)^{-1} \mathcal{S}_t^{-1} = (I - (I+Z_t \mathcal{S}_t)^{-1}Z_t \mathcal{S}_t)\mathcal{S}_t^{-1} \nonumber \\
&= \mathcal{S}_t^{-1} - (I + Z_t \mathcal{S}_t)^{-1} Z_t \succcurlyeq \mathcal{S}_t^{-1} - Z_t, \label{eq:lem_Qconv_step9}
\end{align}
where it was used that $(I+X)^{-1} = I - (I+X)^{-1}X$ and the last inequality again uses Lemma \ref{lem:linalg} (now with $X = Z_t$ and $Y = \mathcal{S}_t$). Therefore, 
\begin{equation}
\mathcal{S}_\infty^{-1}  \succcurlyeq \mathcal{S}_0^{-1} - \sum_{t=1}^\infty Z_t \succcurlyeq (Q_\infty - E_T)^{-1} - W_\infty, \label{eq:lem_Qconv_step10}
\end{equation}
where the latter inequality uses that
\begin{equation}
Z_t = (M_\infty A)^t M_t B R^{-1} B^\top (A^\top M_\infty^\top)^t \preccurlyeq (M_\infty A)^t B R^{-1} B^\top (A^\top M_\infty^\top)^t, \label{eq:lem_Qconv_step11}
\end{equation}
by Lemma \ref{lem:linalg}(iii) (with $X = BR^{-1}B^\top$ and $Y = \mathcal{Q}_t$) and the definition of $W_\infty$ as
\begin{equation}
W_\infty := \sum_{t=0}^\infty (M_\infty A)^t B R^{-1} B^\top (A^\top M_\infty^\top)^t, \label{eq:Grammian_DT} \\
\end{equation}
The result now follows by taking norms in \eqref{eq:lem_Qconv_step6} noting that 
\begin{equation}
\| \mathcal{S}_t \| \leq \| \mathcal{S}_\infty \| \leq \| ((Q_\infty-E_T)^{-1} + W_\infty)^{-1} \|. 
\end{equation}
\end{proof}

\subsection{Stability and convergence}

The convergence result for the DRDE from the previous subsection enables the derivation of a stability condition and convergence results for discrete-time RHC. 

\begin{theorem}[Stability of discrete-time RHC] \label{thm:stabRHC_DT}
There exists a constant $\tilde{K}_1$ independent of $t$, $x_0$, $T$, and $\tau$ such that all $t \in \mathbb{N}$
\begin{equation}
|x^{\mathrm{RHC}}_t| \leq \theta_{T-\tau}^t |x_0|, \label{eq:thm_stabRHC_DT}
\end{equation}
where
\begin{equation}
\theta_{T-\tau} = \| M_\infty A \| + \tilde{K}_1 \| M_\infty A \|^{2(T-\tau)-1}.
\end{equation}
\end{theorem}
Observe that the RHC strategy is stabilizing when $|\theta_{T-\tau}| < 1$. Because $\| M_\infty A \| < 1$, it is easy to see that $|\theta_{T-\tau}| < 1$ for $T - \tau$ sufficiently large.

\begin{proof} 
Using the definitions of $M_\infty$ and $M_t$ in \eqref{eq:def_MinfMt}, \eqref{eq:x_inf_opt_DT}, and \eqref{eq:dyn_xT_DT} can be rewritten as 
\begin{equation}
x^\infty_{t+1} = M_\infty A x^\infty_t, \qquad \qquad x^{\mathrm{RHC}}_{t+1} = M^{T,\tau}_t A x^{\mathrm{RHC}}_t,
\end{equation}
where
\begin{equation}
M^{T,\tau}_t = M_{T-1 - (t\ \mathrm{mod}\ \tau)} = (I + BR^{-1}B^\top\mathcal{Q}_{T-1 - (t\ \mathrm{mod}\ \tau)})^{-1}.
\end{equation}
Note that by definition of $M_\infty$ and $M_t$ in \eqref{eq:def_MinfMt}
\begin{equation}
M_t - M_\infty = M_t (M_\infty^{-1} - M_t^{-1}) M_\infty = M_t B R^{-1}B^\top (Q_\infty - \mathcal{Q}_t) M_\infty,
\end{equation}
so that Lemma \ref{lem:Qconv} shows that
\begin{equation}
\| M_\infty A - M^{T,\tau}_t A \| \leq \tilde{K}_1\| M_\infty A \|^{2(T-\tau) -1}, \label{eq:thm_stabRHC_DT_step3}
\end{equation}
where it was used that $0 \preccurlyeq M_t BR^{-1}B^\top \preccurlyeq BR^{-1}B^\top$ according to Lemma \ref{lem:linalg} and $\tilde{K}_1 = \| BR^{-1}B^\top \|\tilde{K}_0$. The conclusion now follows because
\begin{equation}
\| M^{T,\tau}_t A \| \leq \| M_\infty A \| + \| M^{T,\tau}_t A - M_\infty A \| \leq \theta_{T-\tau}, 
\end{equation}
and $x^{\mathrm{RHC}}_{t+1} = M^{T,\tau}_t A x^{\mathrm{RHC}}_t$. 
%To see that the RHC strategy is stabilizing when $\theta^{T,\tau} < 1$, note that \eqref{eq:dyn_xRH_DT} shows that
%\begin{equation}
%| x^{\mathrm{RHC}}_{t+1}|  = \| A^{T,\tau}_t \| | x^{\mathrm{RHC}}_t |, \quad x^{\mathrm{RHC}}_0 = x^0, \quad \Rightarrow \quad |x^{\mathrm{RHC}}_t| \leq (\theta^{T,\tau})^t |x^0|,
%\end{equation}
%which shows that $x^{\mathrm{RHC}}_t \rightarrow 0$ if $\theta^{T,\tau} < 1$. 
\end{proof}

The following lemma shows that $x^{\mathrm{RHC}}$ converges to $x^\infty$ and that $u^{\mathrm{RHC}}$ converges to $u^*$ for $T - \tau \rightarrow \infty$. 

\begin{theorem}[Convergence of discrete-time RHC] \label{thm:convRHC_DT}
There exist a constant $K$ independent of $t$, $x_0$, $T$, and $\tau$ such that
\begin{equation}
|x^{\mathrm{RHC}}_t - x^\infty_t| + |u^{\mathrm{RHC}}_t - u^\infty_t| \leq K \| M_\infty A \|^{2(T-\tau)+1} (t+1) \theta_{T-\tau}^t |x_0|, \label{eq:thm_convRHCx_DT}
\end{equation}
Furthermore, if $|\theta_{T-\tau}| < 1$ there exists a constant $K$ independent of $t$, $x_0$, $T$, and $\tau$ such that
\begin{equation}
J^\infty(u^{\mathrm{RHC}}) - J^\infty(u^\infty) \leq K  \frac{\| M_\infty A \|^{4(T-\tau)+2}}{(1-\theta_{T-\tau}^2)^3} |x_0|^2.\label{eq:thm_convRHCJ_DT}
\end{equation}
\end{theorem}

\begin{proof} Throughout the proof, $K$ denotes a constant independent of $t$, $x_0$, $T$, and $\tau$ that may vary from line to line. 
Define $e^{\mathrm{RHC}}_t := x^{\mathrm{RHC}}_t - x^\infty_t$, using \eqref{eq:def_MinfMt} it follows that
\begin{equation}
e^{\mathrm{RHC}}_{t+1} = 
M^{T,\tau}_t A x^{\mathrm{RHC}}_t - M_\infty A x^\infty_{t+1} = 
M_\infty A e^{\mathrm{RHC}}_t + (M^{T,\tau}_t A - M_\infty A)x^{\mathrm{RHC}}_t. \label{eq:thm_convRHC_DT_step1}
\end{equation}
Taking norms and using \eqref{eq:thm_stabRHC_DT_step3} shows that
\begin{equation}
|e^{\mathrm{RHC}}_{t+1}| \leq \| M_\infty A \| |e^{\mathrm{RHC}}_t| + K \| M_\infty A \|^{2(T-\tau)+1} |x^{\mathrm{RHC}}_t|.
\label{eq:thm_convRHC_DT_step2}
\end{equation}
By induction over $t$, it is easy to verify that \eqref{eq:thm_convRHC_DT_step2} implies that
\begin{equation}
|e^{\mathrm{RHC}}_t| \leq K \| M_\infty A \|^{2(T-\tau)+1} \sum_{s=0}^{t-1} \| M_\infty A\|^{t-s} |x^{\mathrm{RHC}}_s|
\end{equation}
Inserting the estimate from Theorem \ref{thm:stabRHC_DT} and using that $|\theta_{T-\tau}| > \| M_\infty A \|$, it follows that
\begin{equation}
|e^{\mathrm{RHC}}_t| \leq K \| M_\infty A \|^{2(T-\tau)+1} t \theta_{T-\tau}^t |x_0|. \label{eq:thm_convRHC_DT_step3}
\end{equation}

For the bound on the controls, observe that \eqref{eq:xT_dyn_DT_1} shows that
\begin{equation}
u^{\mathrm{RHC}}_t = -R^{-1}B^\top  \mathcal{Q}_{t'} M_{t'} A x^{\mathrm{RHC}}_t, \qquad u^\infty_t = -R^{-1} B^\top Q_\infty M_\infty A x^\infty_t, 
\end{equation}
where $t' = T-1 - (t \mod \tau)$ for brevity and $M_t$ as in \eqref{eq:def_MinfMt}. 
Therefore,
\begin{align}
u^{\mathrm{RHC}}_t - u^\infty_t &= -R^{-1}B^\top \mathcal{Q}_{t'} M_{t'} A x^{\mathrm{RHC}}_t + R^{-1} B^\top Q_\infty M_\infty A x^\infty_t \nonumber \\
&= -R^{-1}B^\top (Q_\infty M_\infty A e^{\mathrm{RHC}}_t + (\mathcal{Q}_{t'}M_{t'} - Q_\infty M_\infty)A x^{\mathrm{RHC}}_t).
\label{eq:thm_convRHC_DT_step4}
\end{align}
where the definitions of $M_\infty$ and  $M_t$ \eqref{eq:def_MinfMt} have been used. Note that
\begin{align}
\mathcal{Q}_{t'}M_{t'} - Q_\infty M_\infty &= \mathcal{Q}_{t'}M_{t'} (M_\infty^{-1} - M_{t'}^{-1})M_\infty + (\mathcal{Q}_{t'} - Q_\infty)M_\infty \nonumber \\
&= \mathcal{Q}_{t'}M_{t'} BR^{-1}B^\top (Q_\infty - \mathcal{Q}_{t'}) M_\infty + (\mathcal{Q}_{t'} - Q_\infty)M_\infty \nonumber \\
&= (I - \mathcal{Q}_{t'}M_{t'}BR^{-1}B^\top )(\mathcal{Q}_{t'} - Q_\infty)M_\infty \nonumber \\
&= (I - \mathcal{Q}_{t'}BR^{-1}B^\top M_{t'}^\top )(\mathcal{Q}_{t'} - Q_\infty)M_\infty \nonumber \\
&= (I + \mathcal{Q}_{t'}BR^{-1}B^\top )^{-1}(\mathcal{Q}_{t'} - Q_\infty)M_\infty
\label{eq:thm_convRHC_DT_step5}
\end{align}
where the last two identities follows because $M_{t'}BR^{-1}B^\top$ is symmetric by Lemma \ref{lem:linalg} and because $I - X(I+X)^{-1} = (1+X)^{-1}$. Because all eigenvalues of $I + \mathcal{Q}_{t'}BR^{-1}B^\top$ are larger than 1, it follows that
\begin{equation}
\| (\mathcal{Q}_{t'}M_{t'} - Q_\infty M_\infty)Ax^{\mathrm{RHC}}_t \| \leq \| Q_{t'} - Q_\infty \| | M_\infty A x^{\mathrm{RHC}}_t | \label{eq:thm_convRHC_DT_step6}
\end{equation}
Taking norms in \eqref{eq:thm_convRHC_DT_step4}, it follows that
\begin{equation}
|u^{\mathrm{RHC}}_t - u^*_t| \leq K |e^{\mathrm{RHC}}_t| + K \| Q_{t'} - Q_\infty \| | x^{\mathrm{RHC}}_t |.
\label{eq:thm_convRHC_DT_step7}
\end{equation}
The estimate then follows after using \eqref{eq:thm_convRHC_DT_step3} to estimate $|e^{\mathrm{RHC}}_t|$, Lemma \ref{lem:Qconv} to estimate $\| \mathcal{Q}_t - Q_\infty\|$, and Theorem \ref{thm:stabRHC_DT} to estimate $|x^{\mathrm{RHC}}_t|$. 

To prove the second statement, now assume that $|\theta_{T-\tau}| < 1$ and write $v^{\mathrm{RHC}}_t := u^{\mathrm{RHC}}_t - u^\infty_t$. Because $u^\infty$ is the minimizer of $J^\infty$ and $J^\infty$ is quadratic, it follows similarly as in Theorem \ref{thm:RHCconv} that 
\begin{align}
J^\infty(u^{\mathrm{RHC}}) - J^\infty(u^\infty) &= \frac{1}{2} \sum_{t=0}^\infty \left( |Ce^{\mathrm{RHC}}_t|^2 + (v^{\mathrm{RHC}}_t)^\top R v^{\mathrm{RHC}}_t \right) \nonumber \\
&\leq K \sum_{t=0}^\infty |e_t^{\mathrm{RHC}}|^2 + |v^{\mathrm{RHC}}_t|^2. \label{eq:diffJ}
\end{align}
Inserting the estimates from \eqref{eq:thm_convRHCx_DT} and making use of the fact that for $0 \leq \theta < 1$
\begin{equation}
\sum_{t=0}^\infty \theta^{2t} = \frac{1}{1-\theta^2}, \quad \sum_{t=0}^\infty t \theta^{2t-1} = \frac{\theta}{(1-\theta^2)^2}, \quad \sum_{t=0}^\infty t^2 \theta^{2t-2} = \frac{\theta^2+1}{(1-\theta^2)^3},
\end{equation}
\eqref{eq:thm_convRHCJ_DT} follows. 
\end{proof}

\section*{Acknowledgments}
We would like to thank Manuel Schaller for his helpful comment that inspired this paper. 

\bibliographystyle{siamplain}
%\bibliography{references}
\bibliography{MPC_article.bbl}

\end{document}